\documentclass[a4paper,oneside,10pt]{article}%
\usepackage{amsmath}
\usepackage{amsfonts}
\usepackage{amssymb}
\usepackage{graphicx}
\usepackage{color}%
\providecommand{\U}[1]{\protect \rule{.1in}{.1in}}
\pdfoutput=1

\newtheorem{theorem}{Theorem}[section]

\newtheorem{definition}[theorem]{Definition}

\newtheorem{example}[theorem]{Example}

\newtheorem{lemma}[theorem]{Lemma}

\newtheorem{proposition}[theorem]{Proposition}
\newtheorem{remark}[theorem]{Remark}

\newtheorem{sch}{Scheme}
\newenvironment{proof}[1][Proof]{\noindent \textbf{#1.} }{\  \rule{0.5em}{0.5em}}
\numberwithin{equation}{section}
\begin{document}

\title{Discrete-time approximation for backward stochastic differential equations
driven by $G$-Brownian motion}
\author{Lianzi Jiang\thanks{College of Mathematics and Systems Science, Shandong
University of Science and Technology, Qingdao, Shandong 266590, PR China.
jianglianzi95@163.com. Research supported by National Natural Science
Foundation of Shandong Province (No. ZR2023QA090).}
\and Mingshang Hu\thanks{Zhongtai Securities Institute for Financial Studies,
Shandong University, Jinan, Shandong 250100, PR China. humingshang@sdu.edu.cn.
Research supported by National Natural Science Foundation of China (No.
12326603, 11671231) and National Key R\&D Program of China (No.
2018YFA0703900).}}
\date{}
\maketitle

\textbf{Abstract}. In this paper, we study the discrete-time approximation
schemes for a class of backward stochastic differential equations driven by
$G$-Brownian motion ($G$-BSDEs) which corresponds to the hedging pricing of
European contingent claims. By introducing an auxiliary extended
$\widetilde{G}$-expectation space, we propose a class of $\theta$-schemes to
discrete $G$-BSDEs in this space. With the help of nonlinear stochastic
analysis techniques and numerical analysis tools, we prove that our schemes
admit half-order convergence for approximating $G$-BSDE in the general case.
In some special cases, our schemes can achieve a first-order convergence rate.
Finally, we give an implementable numerical scheme for $G$-BSDEs based on
Peng's central limit theorem and illustrate our convergence results with
numerical examples.
%several numerical tests are given to support the theoretical results.



\section{Introduction}

In 1990, Pardoux and Peng in \cite{PP1990} proved the existence and uniqueness
theorem for nonlinear backward stochastic differential equations (BSDEs).
Since then, the theory of BSDEs has grown into an indispensable tool in many
areas. In particular, it is used to measure financial risk under model
uncertainty. But this type of BSDE is difficult to price path-dependent
contingent claims in the uncertain volatility model. To overcome this
obstacle, Peng systemically established a time-consistent $G$-expectation
theory \cite{P2007,P2008,P2010}. Under the $G$-expectation framework, a new
type of Brownian motion called $G$-Brownian motion was constructed and the
corresponding stochastic calculus of It\^{o}'s type was established. Under
standard conditions on $f,g_{ij}$ and $\xi$, Hu et al. \cite{HJPS20141}
obtained the existence and uniqueness of solution $(Y,Z,K)$ for the following
fully nonlinear BSDE driven by $G$-Brownian motion ($G$-BSDE for short, in
this paper we always use Einstein convention):
\begin{equation}
Y_{t}=\xi+\int_{t}^{T}f\left(  s,Y_{s},Z_{s}\right)  ds+\int_{t}^{T}%
g_{ij}\left(  s,Y_{s},Z_{s}\right)  d\langle B^{i},B^{j}\rangle_{s}-\int
_{t}^{T}Z_{s}dB_{s}-\left(  K_{T}-K_{t}\right)  \text{,} \label{G-BSDE_Plimi}%
\end{equation}
where $B_{t}=\left(  B_{t}^{1},\ldots,B_{t}^{d}\right)  ^{\top}$ is a
$d$-dimensional $G$-Brownian motion, $\langle B^{i},B^{j}\rangle$ is the
corresponding mutual variation process, and $K$ is a decreasing $G$%
-martingale. Furthermore, Hu et al. \cite{HJPS20142} studied nonlinear
Feynman-Kac formula under the framework of $G$-expectation, which built the
relationship between $G$-BSDEs and fully nonlinear partial differential
equations (PDEs). In addition, Cheridito et al. \cite{CSTV2007} and Soner et
al. \cite{STZh2012} developed a new type of fully nonlinear FBSDEs, called the
2FBSDEs, which is also associated with a class of fully nonlinear PDEs.

$G$-expectation theory and its applications have been developed in recent
research
\cite{DHP2011,EJ2013,EJ2014,FPSS2019,HJ2013,Krylov2018,PYY2019,Song2017}. In
particular, using the $G$-expectation theory to describe the volatility
uncertainty in the financial market, Vorbrink \cite{V2014} and Hu et al.
\cite{HJ2013} studied the pricing of European contingent claims. They found
that the superhedging and subhedging prices of the given contingent claim
$\xi$ are related to $Y_{t}$, which is the solution of the following $G$-BSDE
\begin{equation}
Y_{t}=\xi-\int_{t}^{T}(r_{s}Y_{s}+a_{s}Z_{s})ds-\int_{t}^{T}b_{s}^{ij}%
Z_{s}d\langle B^{i},B^{j}\rangle_{s}-\int_{t}^{T}Z_{s}dB_{s}-\left(
K_{T}-K_{t}\right)  . \label{BSDE LINEAR}%
\end{equation}
%Indeed, in the Markovian
%case, the superhedging and subhedging prices are the Black-Scholes prices with
%volatilities $\underline{\sigma}^{2}$ and $\bar{\sigma}^{2}$, respectively.
However, the analytic solution $Y_{t}$ of $G$-BSDE (\ref{BSDE LINEAR}) is
often difficult to obtain. Therefore, proposing an numerical approach to
approximate $Y_{t}$ has important practical significance.

There are a lot of works dealing with numerical methods for BSDEs and 2BSDEs,
see, e.g.,
\cite{B1997,BD2007,CR2014,CR2015,CM2014,EHJ2017,FTW2011,GLW2005,GT2016,GZhZh2015,
KZhZh2015,MZ2005,MT2006,RO2016,Z2004,ZhCP2006,ZhWP2009,ZhZhK2017} and
references therein. However, owing to the nonlinear nature of $G$-expectation
framework, many classical numerical analysis techniques are no longer
applicable under this framework. In recent years, some progress has been made
in numerical approximation research under the $G$-expectation framework.
Dolinsky \cite{D2012} studied a discrete-time analog of $G$-expectation and
proved its convergence rate. By solving a class of HJB equations to generate
the $G$-normal distribution, Yang and Zhao \cite{YZh2016} gave a simulation of
$G$-Brownian motion path. Yang and Li \cite{YL2019} studied the exponential
stability of a class of $\theta$-schemes for stochastic differential equations
driven by $G$-Brownian motion ($G$-SDEs). Peng and Zhang \cite{PZ2020}
proposed a Wong-Zakai approximation for Stratonovich type $G$-SDEs pathwisely.
Using the $G$-expectation representation, Hu and Jiang \cite{HJ2021} proposed
a new kind of numerical scheme for solving $G$-BSDEs, and proved the scheme is
at most half-order of convergence. This is applied to the discrete-time
approximation of stochastic optimal control problems under $G$-expectation in
\cite{J2021}. Best we can see, all the above mentioned works are designed with
maximum convergence rate $\frac{1}{2}$.

Motivated by the contingent claim pricing problems, in this work, we devote
ourselves to studying a class of high order discrete approximation schemes for
the following type of $G$-BSDE
\begin{align}
Y_{t}  &  =\xi+\int_{t}^{T}\left(  f\left(  s,Y_{s}\right)  +a_{s}%
Z_{s}\right)  ds+\int_{t}^{T}\left(  g_{ij}\left(  s,Y_{s}\right)  +b_{s}%
^{ij}Z_{s}\right)  d\langle B^{i},B^{j}\rangle_{s}\label{G-BSDE}\\
&  \text{ \  \  \  \ }-\int_{t}^{T}Z_{s}dB_{s}-\left(  K_{T}-K_{t}\right)
,\nonumber
\end{align}
where $a_{s}(\omega),b_{s}^{ij}(\omega):[0,T]\times \Omega \rightarrow
\mathbb{R}^{1\times d}$ are continuous bounded processes in $M_{G}^{2}(0,T)$
(see Definition \ref{Mp}), and $f,g_{ij}=g_{ji}:[0,T]\times \mathbb{R}%
\rightarrow \mathbb{R}$, $i,j=1,\ldots,d$, are deterministic functions.
Different from the idea in \cite{HJ2021}, by constructing an auxiliary
extended $\widetilde{G}$-expectation space, we introduce a forward process and
then propose a new type of $\theta$-scheme for approximating $G$-BSDE
(\ref{G-BSDE}). We point that, this scheme can solve $Y_{t}$ independent of
$Z_{t}$ and $K_{t}$, which provides an efficient numerical method for the
hedging prices of the contingent claim. With the help of nonlinear stochastic
analysis techniques and numerical analysis tools, we rigorously analyze the
errors of the proposed scheme and prove the convergence. It is shown that our
$\theta$-scheme admits a half-order rate of convergence with respect to
$\Delta t$ in the general case. For the case of $\theta_{1}\in \lbrack0,1]$ and
$\theta_{2}=0$, our schemes can reach first-order convergence in the
deterministic case.

%%%%%for the classical BSDE, using the $\theta$-scheme for solving $Y$ achieves a first-order convergence rate for $\theta\in (0,1]$. Specifically, when $\theta=1/2$, the $\theta$-scheme can reach second-order convergence rate for $Y$ under smoothness assumptions, see for example \cite{ZhWP2009}.
%In particular, when the $G$-BSDE is reduced to a BSDE, due to the introduction of the forward process, our $\theta$-scheme solves $Y$ at most with a first-order convergence rate for any $\theta_{i}\in [0,1] (i=1,2)$.

The paper is organized as follows. In Section 2, we recall some basic
notations and results for stochastic calculus under $G$-framework. In Section
3, we propose a class of discrete-time approximation schemes for $G$-BSDEs.
Section 4 is devoted to the convergence analysis of our the schemes. In
Section 5, we give a feasible numerical scheme and provide some numerical
examples to illustrate our theoretical findings.

\section{Preliminaries}

The main purpose of this section is to recall some basic notions and results
of $G$-Expectation, $G$-Brownian motion and $G$-BSDEs, which are needed in the
sequel. More details can refer to \cite{P2007,P2008,P2010} and references therein.

\begin{definition}
Let $\Omega$ be a given set and let $\mathcal{H}$ be a linear space of real
valued functions defined on $\Omega$, satisfies $c\in \mathcal{H}$ for each
constant $c$ and $\left \vert X\right \vert \in \mathcal{H}$ if $X\in \mathcal{H}%
$. $\mathcal{H}$ is considered as the space of random variables. A sublinear
expectation $\mathbb{\hat{E}}$ on $\mathcal{H}$ is a functional $\mathbb{\hat
{E}}$: $\mathcal{H}\rightarrow \mathbb{R}$ satisfying the following properties:
for all $X,Y$ $\in \mathcal{H}$, we have
\begin{enumerate}
\item[$(i)$] Monotonicity: If $X\geq Y$ then $\mathbb{\hat{E}}\left[
X\right]  \geq \mathbb{\hat{E}}\left[  Y\right]  ;$

\item[$(ii)$] Constant preservation: $\mathbb{\hat{E}}\left[  c\right]  =c;$

\item[$(iii)$] Sub-additivity: $\mathbb{\hat{E}}\left[  X+Y\right]
\leq \mathbb{\hat{E}}\left[  X\right]  +\mathbb{\hat{E}}\left[  Y\right]  ;$

\item[$(iv)$] Positive homogeneity: $\mathbb{\hat{E}}\left[  \lambda X\right]
=\lambda \mathbb{\hat{E}}\left[  X\right]  $ for each $\lambda>0.$
\end{enumerate}
\end{definition}
$(\Omega,\mathcal{H},\mathbb{\hat{E})}$ is called a sublinear expectation
space. From the definition of the sublinear expectation $\mathbb{\hat{E}}$,
the following results can be easily obtained.

\begin{proposition}
\label{PropE}For $X,Y$ $\in \mathcal{H}$, we have
\begin{enumerate}
\item[$\left(  i\right)  $] $\mathbb{\hat{E}}\left[  \lambda X\right]
=\lambda^{+}\mathbb{\hat{E}}\left[  X\right]  +\lambda^{-}\mathbb{\hat{E}%
}\left[  -X\right]  ,$ for each $\lambda \in \mathbb{R};$

\item[$\left(  ii\right)  $] If $\mathbb{\hat{E}}\left[  X\right]
=-\mathbb{\hat{E}}\left[  -X\right]  $, then $\mathbb{\hat{E}}\left[
X+Y\right]  =\mathbb{\hat{E}}\left[  X\right]  +\mathbb{\hat{E}}\left[
Y\right]  ;$

\item[$\left(  iii\right)  $] $|\mathbb{\hat{E}}\left[  X\right]
-\mathbb{\hat{E}}\left[  Y\right]  |\leq \mathbb{\hat{E}}\left[  \left \vert
X-Y\right \vert \right]  ,$i.e. $|\mathbb{\hat{E}}\left[  X\right]
-\mathbb{\hat{E}}\left[  Y\right]  |\leq \mathbb{\hat{E}}\left[  X-Y\right]
\vee \mathbb{\hat{E}}\left[  Y-X\right]  ;$

\item[$\left(  iv\right)  $] $\mathbb{\hat{E}}\left[  \left \vert XY\right \vert
\right]  \leq(\mathbb{\hat{E}}\left[  \left \vert X\right \vert ^{p}\right]
)^{1/p}\cdot(\mathbb{\hat{E}}\left[  \left \vert Y\right \vert ^{q}\right]
)^{1/q},$ for $1\leq p,q<\infty$ with $\frac{1}{p}+\frac{1}{q}=1.$
\end{enumerate}
\end{proposition}

\begin{definition}
Let $X_{1}$ and $X_{2}$ be two $n$-dimensional random variables defined
respectively in sublinear expectation spaces $(\Omega_{1},\mathcal{H}%
_{1},\mathbb{\hat{E}}_{1})$ and $(\Omega_{2},\mathcal{H}_{2},\mathbb{\hat{E}%
}_{2})$. They are called identically distributed, denoted by $X_{1}\overset
{d}{=}X_{2}$, if $\mathbb{\hat{E}}_{1}\left[  \varphi(X_{1})\right]
=\mathbb{\hat{E}}_{2}\left[  \varphi(X_{2})\right]  $, for all $\varphi \in
C_{b,Lip}(\mathbb{R}^{n})$, the space of bounded Lipschitz continuous
functions on $\mathbb{R}^{n}$.
\end{definition}

\begin{definition}
In a sublinear expectation space $(\Omega,\mathcal{H},\mathbb{\hat{E})}$, a
$m$-dimensional random variable $Y\ $is said to be independent from another
$n$-dimensional random variable $X$, denoted by $Y\perp \! \! \! \perp X$, if
for every test function $\varphi \in C_{b,Lip}(\mathbb{R}^{n}\times
\mathbb{R}^{m})$ we have $\mathbb{\hat{E}}\left[  \varphi(X,Y)\right]
=\mathbb{\hat{E}[\hat{E}[\varphi(}x,Y\mathbb{)]}_{x=X}]$. Notably, $Y$ is
independent from $X$ does not necessarily imply that $X$ is independent from
$Y$. If $Y\overset{d}{=}X$ and $Y\perp \! \! \! \perp X$, we call $Y$ an
independent copy of $X$.
\end{definition}

Now we recall the definition of $G$-distributed and Peng's central limit theorem.

\begin{definition}
The pair of $\mathbb{R}^{d}\times \mathbb{R}^{k}$-valued random variables
$(\xi,\eta)$ on a sublinear expectation space $(\Omega,\mathcal{H}%
,\mathbb{\hat{E})}$ is called $G$-distributed if for each $a,b\geq0$ we have
\begin{equation}
\left(  a\xi+b\bar{\xi},a\eta+b\bar{\eta}\right)  \overset{d}{=}\left(
\sqrt{a^{2}+b^{2}}\xi,(a+b)\eta \right)  , \label{G}%
\end{equation}
where $(\bar{\xi},\bar{\eta})$ is an independent copy of $(\xi,\eta)$. The
pair $(\xi,\eta)$ associates with the binary function $G:\mathbb{R}^{k}%
\times \mathbb{S}(d)\rightarrow \mathbb{R}$
\begin{equation}
G\left(  p,A\right)  :=\mathbb{\hat{E}}\left[  \frac{1}{2}\left \langle
A\xi,\xi \right \rangle +\left \langle p,\eta \right \rangle \right]
,\ (p,A)\in \mathbb{R}^{k}\times \mathbb{S}(d), \label{G(p,A)}%
\end{equation}
where $\mathbb{S}(d)$ denotes the collection of all $\mathbb{R}^{d\times d}$
symmetric matrices. Specifically, for each $\varphi \in C_{b,Lip}%
(\mathbb{R}^{d\times k})$, the function defined by
\[
u(t,x,y):=\mathbb{\hat{E}}\left[  \varphi(x+\sqrt{t}\xi,y+t\eta)\right]
,\  \ (t,x,y)\in(0,\infty]\times \mathbb{R}^{d}\times \mathbb{R}^{k}%
\]
is the unique viscosity solution of the following parabolic PDE
\[
\partial_{t}u-G(D_{y}u,D_{x}^{2}u)=0,\text{ \ }u|_{t=0}=\varphi.
\]

\end{definition}

\begin{remark}
\label{remark_G} From the representation theorem of sublinear expectations
(cf. Theorem 1.2.1 in \cite{P2010}), it is easy to obtain that there exists a
bounded and closed subset $\Gamma \subset \mathbb{R}^{k}\times \mathbb{S}_{+}(d)$
such that
\[
G(p,A)=\sup_{(q,Q)\in \Gamma}\left[  \frac{1}{2}tr[AQ]+\langle p,q\rangle
\right]  ,\  \text{for\ }(p,A)\in \mathbb{R}^{k}\times \mathbb{S}(d),
\]
where $\mathbb{S}_{+}(d)$ denotes the collection of nonnegative definite
elements in $\mathbb{S}(d)$.
\end{remark}

\begin{remark}
If the pair $(\xi,\eta)$ satisfies (\ref{G}), then
\[
a\xi+b\bar{\xi}\overset{d}{=}\sqrt{a^{2}+b^{2}}\xi,\text{ \ }a\eta+b\bar{\eta
}\overset{d}{=}(a+b)\eta, \text{ \ for }a,b\geq0.
\]
We call $\xi$ the $G$-normally distributed and $\eta$ the maximally
distributed. In this paper, we always assume that $G$ is non-degenerate, i.e.,
there exist some constants $0<\underline{\sigma}^{2}\leq \bar{\sigma}%
^{2}<\infty$ such that $\frac{1}{2}\underline{\sigma}^{2}tr[A-B]$ $\leq$
$\mathbb{\hat{E}}\left[  \frac{1}{2}\left \langle A\xi,\xi \right \rangle
\right]  -\mathbb{\hat{E}}\left[  \frac{1}{2}\left \langle B\xi,\xi
\right \rangle \right]  $ $\leq \frac{1}{2}\bar{\sigma}^{2}tr[A-B]$ for $A\geq
B$. For simplicity, we assume $\sigma^{2}=\underline{\sigma}^{2}\leq
\bar{\sigma}^{2}=1.$
\end{remark}

\begin{theorem}
[\cite{P2010,P2019}]\label{CLT}Let $\left \{  \left(  X_{i},Y_{i}\right)
\right \}  _{i=1}^{\infty}$ be a sequence of $\mathbb{R}^{d}\times
\mathbb{R}^{k}$-valued random variables on a sublinear expectation space
$(\Omega,\mathcal{H},\mathbb{\tilde{E}})$. We assume that $\left(
X_{i+1},Y_{i+1}\right)  \overset{d}{=}\left(  X_{i},Y_{i}\right)  $ and
$\left(  X_{i+1},Y_{i+1}\right)  $ is independent from $\left \{  \left(
X_{1},Y_{1}\right)  ,\ldots,\left(  X_{i},Y_{i}\right)  \right \}  $ for each
$i=1,2,\cdots.$ We further assume that $\mathbb{\tilde{E}}[X_{1}%
]=\mathbb{\tilde{E}}[-X_{1}]=0$ and
\[
\lim \limits_{\gamma \rightarrow+\infty}\mathbb{\tilde{E}}[(|X_{1}|^{2}%
-\gamma)^{+}]=0,\text{ \ }\lim \limits_{\gamma \rightarrow+\infty}%
\mathbb{\tilde{E}}[(|Y_{1}|-\gamma)^{+}]=0.
\]
Then for each function $\varphi \in C(\mathbb{R}^{d}\times \mathbb{R}^{k})$
satisfying linear growth condition, we have
\begin{equation}
\lim_{n\rightarrow \infty}\mathbb{\tilde{E}}\left[  \varphi \left(  \sum
_{i=1}^{n}\frac{X_{i}}{\sqrt{n}},\sum_{i=1}^{n}\frac{Y_{i}}{n}\right)
\right]  =\mathbb{\hat{E}}\left[  \varphi(\xi,\eta)\right]  , \label{8.1}%
\end{equation}
where $\left(  \xi,\eta \right)  $ is a pair of $G$-distributed random
variables under another sublinear expectation $\mathbb{\hat{E}}$ (possibly
different from $\mathbb{\tilde{E}}$), and the corresponding sublinear function
$G:$ $\mathbb{R}^{k}\times \mathbb{S(}d\mathbb{)}$ $\rightarrow \mathbb{R}$ is
defined by
\[
G\left(  p,A\right)  :=\mathbb{\tilde{E}}\left[  \frac{1}{2}\left \langle
AX_{1},X_{1}\right \rangle +\left \langle p,Y_{1}\right \rangle \right]
,\ (p,A)\in \mathbb{R}^{k}\times \mathbb{S(}d\mathbb{)}.
\]

\end{theorem}

%\begin{remark}
%\label{CLTrate}According to \cite{FPSS2019,HL2019}, the above theorem has the
%following convergence rate
%\[
%\left \vert \mathbb{\hat{E}}\left[  \varphi \left(  \sum_{i=1}^{n}\frac{X_{i}
%}{\sqrt{n}},\sum_{i=1}^{n}\frac{Y_{i}}{n}\right)  \right]  -\mathbb{\hat{E}
%}[\varphi(\xi,\eta)]\right \vert \leq Cn^{-\frac{1}{6}},\label{8.3}%
%\]
%where $C$ is a constant depending on $\varphi$, $X_{i}$ and $Y_{i}$.
%\end{remark}

Next, we shall give definitions of $G$-Brownian motion, $G$-expectation and
conditional $G$-expectation.

\begin{definition}
Let $\Omega_{T}=C_{0}([0,T];\mathbb{R}^{d})$, the space of $\mathbb{R}^{d}%
$-valued continuous functions on $[0,T]$ with $\omega_{0}=0$, be endowed with
the supremum norm. Set
\[
Lip(\Omega_{T})=\{ \varphi(B_{t_{1}},\ldots,B_{t_{n}}):n\geq1,0\leq
t_{1}<t_{2}<\cdots<t_{n}\leq T,\varphi \in C_{b,Lip}(\mathbb{R}^{d\times n})\}
\text{.}%
\]
$G$-expectation on $(\Omega_{T},Lip(\Omega_{T}))$ is a sublinear expectation
defined by%
\[
\mathbb{\hat{E}[}X\mathbb{]=\tilde{E}[}\varphi(\sqrt{t_{1}-t_{0}}\xi
_{1},\ldots,\sqrt{t_{m}-t_{m-1}}\xi_{m})\mathbb{]},
\]
for all $X=\varphi(B_{t_{1}}-B_{t_{0}},\ldots,B_{t_{m}}-B_{t_{m-1}})$,
$B_{t}(\omega)=\omega_{t}$, where $\xi_{1},\ldots,\xi_{m}$ are identically
distributed $d$-dimensional $G$-normally distributed random vectors in a
sublinear expectation space $(\tilde{\Omega},\mathcal{\tilde{H}}%
,\mathbb{\tilde{E})}$ such that $\xi_{i+1}$ is independent from $(\xi
_{1},\ldots,\xi_{i})$, $i=1,\ldots,m-1.$ The corresponding canonical process
$B_{t}=(B_{t}^{i})_{i=1}^{d}$ is called a $G$-Brownian motion and $(\Omega
_{T},Lip(\Omega_{T}),\mathbb{\hat{E})}$ is called a $G$-expectation space.
\end{definition}

\begin{definition}
We assume that $\xi \in Lip(\Omega_{T})$ has the representation $\xi=$
$\varphi(B_{t_{1}}-B_{t_{0}},\ldots,B_{t_{m}}-B_{t_{m-1}})$, define the
conditional $G$-expectation $\mathbb{\hat{E}}_{t_{i}}$ of $\xi$, for some
$1\leq i\leq m$%
\[
\mathbb{\hat{E}}_{t_{i}}[\varphi(B_{t_{1}}-B_{t_{0}},\ldots,B_{t_{m}%
}-B_{t_{m-1}})]=\tilde{\varphi}(B_{t_{1}}-B_{t_{0}},\ldots,B_{t_{i}%
}-B_{t_{i-1}}),
\]
where%
\[
\tilde{\varphi}(x_{1},\ldots,x_{i})=\mathbb{\hat{E}}_{t_{i}}[\varphi
(x_{1},\ldots,x_{i},B_{t_{i+1}}-B_{t_{i}},\ldots,B_{t_{m}}-B_{t_{m-1}})].
\]

\end{definition}

Define $\left \Vert \xi \right \Vert _{p,G}=(\mathbb{\hat{E}[}\left \vert
\xi \right \vert ^{p}\mathbb{])}^{1/p},$ for $\xi \in Lip(\Omega_{T})$ and
$p\geq1$, and denote by $L_{G}^{p}(\Omega_{T})$ the completion of
$Lip(\Omega_{T})$ under the norm $\Vert \cdot \Vert_{L_{G}^{p}}$. Then for
$t\in \left[  0,T\right]  $, $\mathbb{\hat{E}}_{t}[\cdot]$ can be extended
continuously to $L_{G}^{1}(\Omega_{T})$. Let $S_{G}^{0}(0,T)=\{h(t,B_{t_{1}%
\wedge t},\ldots,B_{t_{n}\wedge t}):t_{1},\ldots,t_{n}\in \lbrack0,T],h\in
C_{b,Lip}(\mathbb{R\times R}^{d\times n})\}$. For each given $p\geq1$ and
$\eta \in S_{G}^{0}(0,T)$, define $\Vert \eta \Vert_{S_{G}^{p}}=(\mathbb{\hat{E}%
}[\sup_{t\in \lbrack0,T]}|\eta_{t}|^{p}])^{1/p}$ and denote by $S_{G}^{p}(0,T)$
the completion of $S_{G}^{0}(0,T)$ under $\Vert \cdot \Vert_{S_{G}^{p}}$.

%\begin{proposition}
%\label{PropE_t}For $X$ $\in$ $Lip(\Omega_{T})$, we have
%\begin{enumerate}
%\item[(I)] $\mathbb{\hat{E}}_{t}\left[  \xi \right]  =\xi,$ for $\xi \in
%Lip(\Omega_{t});$
%\item[(II)] $\mathbb{\hat{E}}_{t}\left[  \xi X\right]  =\xi^{+}\mathbb{\hat
%{E}}_{t}\left[  X\right]  +\xi^{-}\mathbb{\hat{E}}_{t}\left[  -X\right]  ,$
%for $\xi \in Lip(\Omega_{t});$
%\item[(III)] $\mathbb{\hat{E}}_{s}[\mathbb{\hat{E}}_{t}\left[  X\right]
%]=\mathbb{\hat{E}}_{s}\left[  X\right]  ,$ for $s\leq t;$
%\item[(IV)] $\mathbb{\hat{E}[\hat{E}}_{t}\left[  X\right]  ]=\mathbb{\hat{E}%
%}\left[  X\right]  .$
%\end{enumerate}
%\end{proposition}

\begin{definition}
\label{Mp}Let $M_{G}^{0}(0,T)$ be the collection of processes in the following
form: for a given partition $\{t_{0},t_{1},\ldots,t_{N}\}$ of $[0,T]$,
\[
\eta_{t}(\omega)=\sum_{n=0}^{N-1}\xi_{n}(\omega)I_{[t_{n},t_{n+1})}(t),
\]
where $\xi_{n}\in Lip(\Omega_{t_{n}})$, $n=0,1,\ldots,N-1$. For each given
$p\geq1$, define $\Vert \eta \Vert_{M_{G}^{p}}=(\mathbb{\hat{E}}[\int_{0}
^{T} \vert \eta_{s} \vert ^{p}ds])^{1/p}$ for $\eta \in M_{G}%
^{0}(0,T)$, and denote by $M_{G}^{p}(0,T)$ the completion of $M_{G}^{0}(0,T)$
under $\Vert \cdot \Vert_{M_{G}^{p}}$.
\end{definition}

For every $\xi_{t}\in$ $M_{G}^{2}(0,T;\mathbb{R}^{d})$, we define
\[
\int_{0}^{T}\xi_{t}dB_{t}:=\sum_{i=1}^{d}\int_{0}^{T}\xi_{t}^{i}dB_{t}%
^{i}\text{.}%
\]
For every $i,j=1,\ldots,d$, the mutual variation of $B^{i}$ and $B^{j}$
\[
\langle B^{i},B^{j}\rangle_{t}:=B_{t}^{i}B_{t}^{j}-\int_{0}^{t}B_{s}^{i}%
dB_{s}^{j}-\int_{0}^{t}B_{s}^{j}dB_{s}^{i}%
\]
is also defined since $B^{i},B^{j}\in M_{G}^{2}(0,T)$. Denote by $\langle
B\rangle_{t}:=(\langle B^{i},B^{j}\rangle_{t})_{i,j=1}^{d}$. Then for
$\eta_{t}\in$ $M_{G}^{1}(0,T;\mathbb{R}^{d})$, the $G$-It\^{o} integral
\[
\int_{0}^{T}\eta_{t}d\langle B\rangle_{t}:=\sum_{i,j=1}^{d}\int_{0}^{T}%
\eta_{t}^{ij}d\langle B^{i},B^{j}\rangle_{t}%
\]
is well defined. For more details of the $G$-It\^{o} integral, one may refer
to \cite{P2008,P2010}.

In this paper, we suppose that the terminal value $\xi$ of\ the $G$-BSDE
(\ref{G-BSDE}) is $\varphi(B_{T})$ and the coefficients satisfy the following conditions:

\begin{enumerate}
\item[(H1)] $f$, $g_{ij}$ are continuous in $t$;

\item[(H2)] There exists a constant $L>0$ such that%
\[%
\begin{array}
[c]{r}%
|\varphi(y)-\varphi(y^{\prime})|+|f\left(  s,y\right)  -f\left(  s,y^{\prime
}\right)  |+\sum \limits_{i,j=1}^{d}|g_{ij}\left(  s,y\right)  -g_{ij}\left(
s,y^{\prime}\right)  |\leq L|y-y^{\prime}|.
\end{array}
\]

\end{enumerate}

From Theorem 4.1 in \cite{HJPS20141}, the $G$-BSDE (\ref{G-BSDE}) has a unique
solution $(Y,Z,K)\in \mathfrak{S}_{G}^{2}(0,T)$, where $\mathfrak{S}_{G}%
^{2}(0,T)$ is the collection of processes $(Y,Z,K)$ such that $Y\in S_{G}%
^{2}(0,T),Z\in M_{G}^{2}(0,T)$, and $K$ is a non-increasing $G$-martingale with
$K_{0}=0$ and $K_{T}\in L_{G}^{2}(\Omega_{T})$.

Now we present the nonlinear Feynman-Kac formula for the $G$-BSDE
(\ref{G-BSDE}).

\begin{theorem}
[\cite{HJPS20142}]\label{Nonlinear Feynman-Kac}Assume that the above
conditions hold. Then the unique solution of $G$-BSDE (\ref{G-BSDE}) can be
represented as $Y_{t}=u\left(  t,B_{t}\right)  $, where $u\left(  t,x\right)
\in C^{\frac{1}{2},1}$ is the unique viscosity solution of the following PDE:
\begin{equation}
\left \{
\begin{array}
[c]{l}%
\partial_{t}u+F\left(  D_{x}^{2}u,D_{x}u,u,t\right)  =0,\\
u\left(  T,x\right)  =\varphi \left(  x\right)  ,
\end{array}
\right.  \label{F-K}%
\end{equation}
where%
\begin{align*}
F\left(  D_{x}^{2}u,D_{x}u,u,t\right)   &  =G\left(  H\left(  D_{x}^{2}%
u,D_{x}u,u,t\right)  \right)  +f\left(  t,u\right)  +a_{t}D_{x}u,\\
H_{ij}\left(  D_{x}^{2}u,D_{x}u,u,t\right)   &  =\left(  D_{x}^{2}u\right)
_{ij}+2\left(  g_{ij}\left(  t,u\right)  +b_{t}^{ij}D_{x}u\right)  .
\end{align*}

\end{theorem}

For readers' convenience, we end this section with our main notations of this paper.

\begin{itemize}
\item $\left \vert \cdot \right \vert :$ the standard Euclidean norm in
$\mathbb{R}$ or $\mathbb{R}^{d}$.

\item $C^{\frac{1}{2},1}:$ continuous functions $\phi:\left[  0,T\right]
\times \mathbb{R}^{d}\rightarrow \mathbb{R}$ such that for any $t,s\in [0,T]$,
$x,y\in \mathbb{R}^{d}$, $|\phi(t,x)-\phi(s,y)|\leq C(|t-s|^{\frac{1}{2}}+|x-y|)$ with $C>0$. $($Note: $\phi(x)\in C^{\frac{1}{2},1}$ means
that $\phi$ is a Lipschitz continuous function$).$

\item $C_{b}^{l,k}:$ continuously differentiable functions $\phi:\left[
0,T\right]  \times$ $\mathbb{R}^{d}\rightarrow \mathbb{R}$ with uniformly
bounded partial derivatives $\partial_{t}^{l_{1}}\phi$ and $\partial
_{x}^{k_{1}}\phi$ for $l_{1}\leq l$ and $k_{1}\leq k$.

\item $\mathbb{\hat{E}}_{t_{n}}^{x}[X]$ $:$ the conditional $\tilde{G}%
$-expectation of $X$, under the $\sigma$-field $\sigma \{B_{t_{n}}=x\}$, i.e.
$\mathbb{\hat{E}}\left[  X|B_{t_{n}}=x\right]  $.
\end{itemize}

\section{Discrete-time approximation for $G$-BSDEs}

\subsection{Explicit expression of $Y_{t}$}

We first consider the case of $d=1$. The results for the multi-dimensional
case will be given later. In order to get the explicit expression of $Y_{t}$
in \eqref{G-BSDE} excluding $Z_{t}$ and $K_{t}$, we construct an auxiliary
extended $\tilde{G}$-expectation space $\,(\tilde{\Omega}_{T},L_{\tilde{G}%
}^{1}(\tilde{\Omega}_{T}),\mathbb{\hat{E}}^{\tilde{G}})$ with $\tilde{\Omega
}_{T}=C_{0}([0,\infty),\mathbb{R}^{2})$ and%
\[
\tilde{G}\left(  A\right)  =\frac{1}{2}\sup_{\sigma^{2}\leq v\leq1}tr\left[
A\left[
\begin{array}
[c]{rr}%
v & 1\\
1 & v^{-1}%
\end{array}
\right]  \right]  ,\text{ \ }A\in \mathbb{S(}2\mathbb{)}\text{.}%
\]
To simplify presentation, we still denote by $\mathbb{\hat{E}}$ the extended
$\tilde{G}$-expectation in the sequel. Let $(B_{t},\tilde{B}_{t})_{t\geq0}$ be
the canonical process in the extended space, it is easy to check that $\langle
B,\tilde{B}\rangle_{t}=t.$

Introduce the following linear $\tilde{G}$-SDE:
\begin{equation}
X_{t}=1+\int_{0}^{t}a_{s}X_{s}d\tilde{B}_{s}+\int_{0}^{t}b_{s}X_{s}%
dB_{s},\text{ \ }t\in \lbrack0,T]. \label{G-SDE}%
\end{equation}
It is easy to verify that
\begin{equation}
X_{t}=\exp \left(  \int_{0}^{t}-a_{s}b_{s}ds\right)  \exp \left(  \int_{0}
^{t}b_{s}dB_{s}-\frac{1}{2}\int_{0}^{t}b_{s}^{2}d\langle B\rangle_{s}\right)
\exp \left(  \int_{0}^{t}a_{s}d\tilde{B}_{s}-\frac{1}{2}\int_{0}^{t}a_{s}%
^{2}d\langle \tilde{B}\rangle_{s}\right)  .
\end{equation}
Applying $G$-It\^{o}'s formula to $X_{t}Y_{t}$, we obtain%
\begin{align*}
X_{t}Y_{t}  &  =X_{T}\xi+\int_{t}^{T}X_{s}f\left(  s,Y_{s}\right)  ds+\int
_{t}^{T}X_{s}g\left(  s,Y_{s}\right)  d\langle B\rangle_{s}-\int_{t}^{T}%
a_{s}X_{s}Y_{s}d\tilde{B}_{s}\\
&  \text{ \  \ }-\int_{t}^{T}(X_{s}Z_{s}+b_{s}X_{s}Y_{s})dB_{s}-\int_{t}%
^{T}X_{s}dK_{s}\text{.}%
\end{align*}
From Lemma 3.4 in \cite{HJPS20141}, we know that $\{ \int_{0}^{t}X_{s}%
dK_{s}\}_{0\leq t\leq T}$ is a $\tilde{G}$-martingale. Thus we get the
explicit expression of $Y_{t}$
\begin{equation}
Y_{t}=\mathbb{\hat{E}}_{t}\left[  X_{T}^{t}\xi+\int_{t}^{T}X_{s}^{t}f\left(
s,Y_{s}\right)  ds+\int_{t}^{T}X_{s}^{t}g\left(  s,Y_{s}\right)  d\langle
B\rangle_{s}\right]  \text{,} \label{2.1}%
\end{equation}
where $X_{s}^{t}=X_{s}/X_{t}$.

\subsection{The discrete-time scheme}

For the time interval $\left[  0,T\right]  $, we introduce the following
partition:%
\begin{equation}
0=t_{0}<t_{1}<\cdots<t_{N}=T \label{time partition}%
\end{equation}
with $\Delta t_{n}=t_{n+1}-t_{n}\ $and $\Delta t=\max \limits_{0\leq n\leq
N-1}\Delta t_{n}$. We denote $\Delta B_{n+1}=B_{t_{n+1}}-B_{t_{n}}\sim N(\{0\}
\times \Delta t_{n}\Sigma)$ and $\Delta \tilde{B}_{n+1}=\tilde{B}_{t_{n+1}%
}-\tilde{B}_{t_{n}}\sim N(\{0\} \times \Delta t_{n}\tilde{\Sigma})$ with
$\Sigma=[\sigma^{2},1]$ and $\tilde{\Sigma}=[1,\frac{1}{\sigma^{2}}]$,
respectively, and $\Delta \langle B\rangle_{n+1}=\langle B\rangle_{t_{n+1}%
}-\langle B\rangle_{t_{n}}$ and $\Delta \langle \tilde{B}\rangle_{n+1}%
=\langle \tilde{B}\rangle_{t_{n+1}}-\langle \tilde{B}\rangle_{t_{n}}$.

Let $f_{s}=f\left(  s,Y_{s}\right)  $ and $g_{s}=g\left(  s,Y_{s}\right)  $.
From \eqref{2.1}, we can obtain
\begin{equation}
Y_{t_{n}}=\mathbb{\hat{E}}_{t_{n}}^{x}\left[  X_{t_{n+1}}^{t_{n}}Y_{t_{n+1}%
}+\int_{t_{n}}^{t_{n+1}}X_{s}^{t_{n}}f\left(  s,Y_{s}\right)  ds+\int_{t_{n}%
}^{t_{n+1}}X_{s}^{t_{n}}g\left(  s,Y_{s}\right)  d\langle B\rangle_{s}\right]
, \label{3.2.1}%
\end{equation}
where%
\begin{align}
X_{t}^{t_{n}}  &  =\exp \left(  -\int_{t_{n}}^{t}a_{s}b_{s}ds\right)
\exp \left(  \int_{t_{n}}^{t}b_{s}dB_{s}-\frac{1}{2}\int_{t_{n}}^{t}b_{s}%
^{2}d\langle B\rangle_{s}\right) \label{3.2.6}\\
&  \text{ \  \ }\exp \left(  \int_{t_{n}}^{t}a_{s}d\tilde{B}_{s}-\frac{1}{2}%
\int_{t_{n}}^{t}a_{s}^{2}d\langle \tilde{B}\rangle_{s}\right)  \text{.}%
\nonumber
\end{align}
It is easy to verify that $X_{t}^{t_{n}}$ satisfies the following $\tilde{G}%
$-SDE
\begin{equation}
X_{t}^{t_{n}}=1+\int_{t_{n}}^{t}a_{s}X_{s}^{t_{n}}d\tilde{B}_{s}+\int_{t_{n}%
}^{t}b_{s}X_{s}^{t_{n}}dB_{s},\text{ \ }t\in \lbrack t_{n},t_{n+1}].
\label{3.2.5}%
\end{equation}
Using the $\theta$-method to approximate the stochastic integrals in
(\ref{3.2.1}), we have
\begin{align}
Y_{t_{n}}  &  =\mathbb{\hat{E}}_{t_{n}}^{x}\left[  X_{t_{n+1}}^{t_{n}%
}Y_{t_{n+1}}+\left(  \theta_{1}f_{t_{n}}+(1-\theta_{1})X_{t_{n+1}}^{t_{n}%
}f_{t_{n+1}}\right)  \Delta t_{n}\right. \nonumber \\
&  \text{ \  \ }+\left.  \left(  \theta_{2}g_{t_{n}}+(1-\theta_{2})X_{t_{n+1}%
}^{t_{n}}g_{t_{n+1}}\right)  \Delta \langle B\rangle_{n+1}\right]  +R_{y}^{n},
\label{3.2.3}%
\end{align}
where the deterministic parameters $\theta_{i}\in \left[  0,1\right]  (i=1,2)$,
and $R_{y}^{n}=R_{y_{1}}^{n}+R_{y_{2}}^{n}$ with
\begin{align}
R_{y_{1}}^{n}  &  =\mathbb{\hat{E}}_{t_{n}}^{x}\left[  X_{t_{n+1}}^{t_{n}%
}Y_{t_{n+1}}+\int_{t_{n}}^{t_{n+1}}X_{s}^{t_{n}}f_{s}ds+\int_{t_{n}}^{t_{n+1}%
}X_{s}^{t_{n}}g_{s}d\langle B\rangle_{s}\right] \label{RY1}\\
&  \text{ \  \ }-\mathbb{\hat{E}}_{t_{n}}^{x}\left[  X_{t_{n+1}}^{t_{n}%
}Y_{t_{n+1}}+\left(  \theta_{1}f_{t_{n}}+(1-\theta_{1})X_{t_{n+1}}^{t_{n}%
}f_{t_{n+1}}\right)  \Delta t_{n}+\int_{t_{n}}^{t_{n+1}}X_{s}^{t_{n}}%
g_{s}d\langle B\rangle_{s}\right] \nonumber
\end{align}
and%
\begin{align}
R_{y_{2}}^{n}  &  =\mathbb{\hat{E}}_{t_{n}}^{x}\left[  X_{t_{n+1}}^{t_{n}%
}Y_{t_{n+1}}+\left(  \theta_{1}f_{t_{n}}+(1-\theta_{1})X_{t_{n+1}}^{t_{n}%
}f_{t_{n+1}}\right)  \Delta t_{n}+\int_{t_{n}}^{t_{n+1}}X_{s}^{t_{n}}%
g_{s}d\langle B\rangle_{s}\right] \nonumber \\
&  \text{ \  \ }-\mathbb{\hat{E}}_{t_{n}}^{x}\left[  X_{t_{n+1}}^{t_{n}%
}Y_{t_{n+1}}+\left(  \theta_{1}f_{t_{n}}+(1-\theta_{1})X_{t_{n+1}}^{t_{n}%
}f_{t_{n+1}}\right)  \Delta t_{n}\right. \label{RY2}\\
&  \text{ \  \ }+\left.  \left(  \theta_{2}g_{t_{n}}+(1-\theta_{2})X_{t_{n+1}%
}^{t_{n}}g_{t_{n+1}}\right)  \Delta \langle B\rangle_{n+1}\right]  .\nonumber
\end{align}
In the meanwhile, we employ Euler's method to approximate the stochastic
integrals in $X_{t_{n+1}}^{t_{n}}$, denoted by $\tilde{X}_{t_{n+1}}^{t_{n}}$,
then
\begin{align}
\tilde{X}_{t_{n+1}}^{t_{n}}  &  =\exp \left(  -\int_{t_{n}}^{t_{n+1}}a_{t_{n}%
}b_{t_{n}}ds\right)  \exp \left(  \int_{t_{n}}^{t_{n+1}}b_{t_{n}}dB_{s}%
-\frac{1}{2}\int_{t_{n}}^{t_{n+1}}b_{t_{n}}^{2}d\langle B\rangle_{s}\right)
\label{3.2.4}\\
&  \text{ \  \ }\exp \left(  \int_{t_{n}}^{t_{n+1}}a_{t_{n}}d\tilde{B}_{s}%
-\frac{1}{2}\int_{t_{n}}^{t_{n+1}}a_{t_{n}}^{2}d\langle \tilde{B}\rangle
_{s}\right)  .\nonumber
\end{align}
Replacing $X_{t_{n+1}}^{t_{n}}$ with $\tilde{X}_{t_{n+1}}^{t_{n}}$ in
(\ref{3.2.3}), we get the following reference equation%
\begin{align}
Y_{t_{n}}  &  =\mathbb{\hat{E}}_{t_{n}}^{x}\left[  \tilde{X}_{t_{n+1}}^{t_{n}%
}Y_{t_{n+1}}+\left(  \theta_{1}f_{t_{n}}+(1-\theta_{1})\tilde{X}_{t_{n+1}%
}^{t_{n}}f_{t_{n+1}}\right)  \Delta t_{n}\right. \nonumber \\
&  \text{ \  \ }+\left.  \left(  \theta_{2}g_{t_{n}}+(1-\theta_{2})\tilde
{X}_{t_{n+1}}^{t_{n}}g_{t_{n+1}}\right)  \Delta \langle B\rangle_{n+1}\right]
+R_{x}^{n}+R_{y}^{n}, \label{reference}%
\end{align}
where%
\begin{align}
R_{x}^{n}  &  =\mathbb{\hat{E}}_{t_{n}}^{x}\left[  X_{t_{n+1}}^{t_{n}%
}Y_{t_{n+1}}+\left(  \theta_{1}f_{t_{n}}+(1-\theta_{1})X_{t_{n+1}}^{t_{n}%
}f_{t_{n+1}}\right)  \Delta t_{n}\right. \nonumber \\
&  \text{ \  \ }+\left.  \left(  \theta_{2}g_{t_{n}}+(1-\theta_{2})X_{t_{n+1}%
}^{t_{n}}g_{t_{n+1}}\right)  \Delta \langle B\rangle_{n+1}\right] \label{RX}\\
&  \text{ \  \ }-\mathbb{\hat{E}}_{t_{n}}^{x}\left[  \tilde{X}_{t_{n+1}}%
^{t_{n}}Y_{t_{n+1}}+\left(  \theta_{1}f_{t_{n}}+(1-\theta_{1})\tilde
{X}_{t_{n+1}}^{t_{n}}f_{t_{n+1}}\right)  \Delta t_{n}\right. \nonumber \\
&  \text{ \  \ }+\left.  \left(  \theta_{2}g_{t_{n}}+(1-\theta_{2})\tilde
{X}_{t_{n+1}}^{t_{n}}g_{t_{n+1}}\right)  \Delta \langle B\rangle_{n+1}\right]
.\nonumber
\end{align}

Let $Y^{n}$ denote the discrete-time approximation to $Y_{t}$ at time level
$t=t_{n}$, and use $X^{n}$ to denote $\tilde{X}_{t_{n}}^{t_{n-1}}$, for
$n=0,1,\ldots,N$. Based on (\ref{reference}), we propose the following
$\theta$-scheme for solving $Y_{t}$ of $G$-BSDE (\ref{G-BSDE}), which is
independent of $Z_{t}$ and $K_{t}$.

\begin{sch}
\label{Scheme1} Given random variables $Y^{N}$, solve random variables $Y^{n}%
$, $n=N-1,\ldots,1,0$ from%
\begin{equation}
\left \{
\begin{array}
[c]{l}%
Y^{n}=\mathbb{\hat{E}}_{t_{n}}^{x}\left[  X^{n+1}Y^{n+1}+\left(  \theta
_{1}f(t_{n},Y^{n})+(1-\theta_{1})X^{n+1}f(t_{n+1},Y^{n+1})\right)  \Delta
t_{n}\right. \\
\text{ \  \  \  \  \ }+\left.  \left(  \theta_{2}g(t_{n},Y^{n})+(1-\theta
_{2})X^{n+1}g(t_{n+1},Y^{n+1})\right)  \Delta \langle B\rangle_{n+1}\right]
,\\
X^{n+1}=\exp \left(  -a_{t_{n}}b_{t_{n}}\Delta t_{n}\right)  \exp(b_{t_{n}%
}\Delta B_{n+1}-\frac{1}{2}b_{t_{n}}^{2}\Delta \langle B\rangle_{n+1})\\
\text{ \  \  \  \ }\  \exp(a_{t_{n}}\Delta \tilde{B}_{n+1}-\frac{1}{2}a_{t_{n}}%
^{2}\Delta \langle \tilde{B}\rangle_{n+1}),
\end{array}
\right.  \label{scheme1}%
\end{equation}
with the deterministic parameters $\theta_{i}\in \left[  0,1\right]  (i=1,2)$.
\end{sch}

%\begin{remark}
%In general, we may pay more attention to the solution $Y_{t}$ of the $G$-BSDE
%(\ref{G-BSDE}), which reflects the hedging price of a given contingent claim.
%By introducing the forward process $X$, Scheme \ref{Scheme1} can achieve the
%approximation of solution $Y_{t}$ independent of $Z_{t}$ and $K_{t}$.
%\end{remark}

\begin{remark}
By choosing different $\theta_{i}(i=1,2)\in \lbrack0,1]$, we obtain a class of
different $\theta$-schemes for $G$-BSDE. In particular, our $\theta$-schemes
can admit a first order rate of convergence when $\theta_{1}\in \lbrack0,1]$
and $\theta_{2}=0$ in the deterministic case, which will be proved in the
following section.
\end{remark}

\subsection{\textcolor[rgb]{1.00,0.00,0.00}{Multi-dimensional $G$-Brownian motion case}}%

In this subsection, we shall extend the discrete-time approximation scheme to
the multi-dimensional case. We first take the two-dimensional case as an example.

\begin{example}
Let $G:$ $\mathbb{S}(2)\rightarrow \mathbb{R}$ be a given sublinear function
such that%
\[
G\left(  A\right)  =\frac{1}{2}\sup_{Q\in \Sigma}tr[AQ],
\]
with $\Sigma \subset$ $\mathbb{S}_{+}(2)$ being a bound and closed set of
diagonal matrices, and let $B_{t}=(B_{t}^{1},B_{t}^{2})^{\top}$ be the
$2$-dimensional $G$-Brownian motion on the corresponding $G$-expectation space
$(\Omega_{T},L_{G}^{1}(\Omega_{T}),\mathbb{\hat{E}})$ with $\Omega
_{T}=C_{0}([0,T];\mathbb{R}^{2})$. From \cite[Corollary 3.5.8]{P2010}, we know
that $\langle B^{i},B^{j}\rangle=0$ for $i\neq j$. In this case, $G$-BSDE
(\ref{G-BSDE}) becomes to
\begin{align*}
Y_{t}  &  =\xi+\int_{t}^{T}\left(  f\left(  s,Y_{s}\right)  +a_{s}%
Z_{s}\right)  ds+\int_{t}^{T}\left(  g_{ii}\left(  s,Y_{s}\right)  +b_{s}%
^{i}Z_{s}\right)  d\langle B^{i},B^{i}\rangle_{s}\\
&  \text{ \  \ }-\int_{t}^{T}Z_{s}dB_{s}-\left(  K_{T}-K_{t}\right)  ,
\end{align*}
with $a_{s}=(a_{s}^{1},a_{s}^{2})$, $b_{s}^{i}=(b_{s}^{i1},b_{s}^{i2})$,
$i=1,2$.

\textcolor[rgb]{1.00,0.00,0.00}{For any $Q=\bigg(\begin{array}
[c]{cc}r_{1} & 0\\
0 & r_{2}\end{array}
\bigg)\in \Sigma$ with $r_{i}>0$, $i=1,2$, by choosing $s_{1}=\sup_{Q\in \Sigma
}\left(  \frac{r_{1}^{2}}{r_{2}}+\frac{r_{1}^{2}}{r_{2}^{2}}\right)  +1$ \ and
\ $s_{2}=\sup_{Q\in \Sigma}\left(  \frac{r_{2}^{2}}{r_1}+\frac{r_{2}^{2}}{r_{1}^{2}}\right)  +1$, we can define the following nonnegative definite
matrix}
\[
\tilde{Q}=\left(
\begin{array}
[c]{cccccc}%
r_{1} & 0 & 1 & 0 & 0 & r_{2}\\
0 & r_{2} & 0 & 1 & r_{1} & 0\\
1 & 0 & r_{1}^{-1}+1 & 0 & 0 & 0\\
0 & 1 & 0 & r_{2}^{-1}+1 & 0 & 0\\
0 & r_{1} & 0 & 0 & s_{1} & 0\\
r_{2} & 0 & 0 & 0 & 0 & s_{2}%
\end{array}
\right)  .
\]
Then we construct an auxiliary extended $\tilde{G}$-expectation space
$\,(\tilde{\Omega}_{T},L_{\tilde{G}}^{1}(\tilde{\Omega}_{T}),\mathbb{\hat{E}%
})$ with $\tilde{\Omega}_{T}=C_{0}([0,T],\mathbb{R}^{6})$ and%
\[
\tilde{G}\left(  A\right)  =\frac{1}{2}\sup_{Q\in \Sigma}tr\big[A\tilde
{Q}\big],\text{ for }A\in \mathbb{S}(6).
\]
Let $(B_{t}^{1},B_{t}^{2},\tilde{B}_{t}^{1},\tilde{B}_{t}^{2},\tilde{B}%
_{t}^{3},\tilde{B}_{t}^{4})^{\top}$ be the canonical process in the extended
space. It is easy to check that $\langle \tilde{B}^{1},B^{1}\rangle_{t}%
=\langle \tilde{B}^{2},B^{2}\rangle_{t}=t$, $\langle \tilde{B}^{1},B^{2}%
\rangle_{t}=\langle \tilde{B}^{2},B^{1}\rangle_{t}=0$, $\langle \tilde{B}%
^{3},B^{1}\rangle_{t}=\langle \tilde{B}^{4},B^{2}\rangle_{t}=0$, $\langle
\tilde{B}^{3},B^{2}\rangle_{t}=\langle B^{1},B^{1}\rangle_{t}$, and
$\langle \tilde{B}^{4},B^{1}\rangle_{t}=\langle B^{2},B^{2}\rangle_{t}.$

Introduce the following linear $\tilde{G}$-SDE on $[0,T]$
\[
X_{t}=1+\int_{0}^{t}b_{s}^{11}X_{s}dB_{s}^{1}+\int_{0}^{t}b_{s}^{22}%
X_{s}dB_{s}^{2}+\int_{0}^{t}a_{s}^{1}X_{s}d\tilde{B}_{s}^{1}+\int_{0}^{t}%
a_{s}^{2}X_{s}d\tilde{B}_{s}^{2}+\int_{0}^{t}b_{s}^{12}X_{s}d\tilde{B}_{s}%
^{3}+\int_{0}^{t}b_{s}^{21}X_{s}d\tilde{B}_{s}^{4}.
\]
It is easy to verify that%
\begin{align*}
&  X_{t}=\exp \left(  \int_{0}^{t}b_{s}^{11}dB_{s}^{1}-\frac{1}{2}\int_{0}%
^{t}|b_{s}^{11}|^{2}d\langle B^{1}\rangle_{s}\right)  \exp \left(  \int_{0}%
^{t}b_{s}^{22}dB_{s}^{2}-\frac{1}{2}\int_{0}^{t}|b_{s}^{22}|^{2}d\langle
B^{2}\rangle_{s}\right) \\
&  \exp \left(  \int_{0}^{t}a_{s}^{1}d\tilde{B}_{s}^{1}-\frac{1}{2}\int_{0}%
^{t}|a_{s}^{1}|^{2}d\langle \tilde{B}^{1}\rangle_{s}\right)  \exp \left(
\int_{0}^{t}a_{s}^{2}d\tilde{B}_{s}^{2}-\frac{1}{2}\int_{0}^{t}|a_{s}^{2}%
|^{2}d\langle \tilde{B}^{2}\rangle_{s}\right) \\
&  \exp \left(  \int_{0}^{t}b_{s}^{12}d\tilde{B}_{s}^{3}-\frac{1}{2}\int
_{0}^{t}|b_{s}^{12}|^{2}d\langle \tilde{B}^{3}\rangle_{s}\right)  \exp \left(
\int_{0}^{t}b_{s}^{21}d\tilde{B}_{s}^{4}-\frac{1}{2}\int_{0}^{t}|b_{s}%
^{21}|^{2}d\langle \tilde{B}^{4}\rangle_{s}\right) \\
&  \exp \left(  -\int_{0}^{t}b_{s}^{12}b_{s}^{22}d\langle B^{1}\rangle_{s}%
-\int_{0}^{t}b_{s}^{11}b_{s}^{21}d\langle B^{2}\rangle_{s}\right)  \exp \left(
-\int_{0}^{t}(a_{s}^{1}b_{s}^{11}+a_{s}^{2}b_{s}^{22})\text{d}s\right)  ,
\end{align*}
where $\langle B^{i}\rangle_{t}:=\langle B^{i},B^{i}\rangle_{t}$ and
$\langle \tilde{B}^{i}\rangle_{t}:=\langle \tilde{B}^{i},\tilde{B}^{i}%
\rangle_{t}$, for $i=1,2,3,4$. Applying $G$-It\^{o}'s formula to $X_{t}Y_{t}$, we
obtain%
\begin{align*}
X_{t}Y_{t}  &  =X_{T}\xi+\int_{t}^{T}X_{s}f\left(  s,Y_{s}\right)
\text{d}s+\int_{t}^{T}X_{s}g_{ii}\left(  s,Y_{s}\right)  d\langle B^{i}%
\rangle_{s}-\int_{t}^{T}X_{s}Z_{s}dB_{s}\\
&  \text{ \  \ }-\int_{t}^{T}Y_{s}dX_{s}-\int_{t}^{T}X_{s}dK_{s},
\end{align*}
which implies
\[
Y_{t}=\mathbb{\hat{E}}_{t}\left[  X_{T}^{t}\xi+\int_{t}^{T}X_{s}^{t}f\left(
s,Y_{s}\right)  \text{d}s+\int_{t}^{T}X_{s}^{t}g_{ii}\left(  s,Y_{s}\right)
d\langle B^{i}\rangle_{s}\right]  \text{,}%
\]
where $X_{s}^{t}=X_{s}/X_{t}$. Similar to the one-dimensional case, using
$\theta$-method, we can get the following reference equation%

\begin{align}
Y_{t_{n}}  &  =\mathbb{\hat{E}}_{t_{n}}^{x}\left[  \tilde{X}_{t_{n+1}}^{t_{n}%
}Y_{t_{n+1}}+\left(  \theta_{1}f_{t_{n}}+(1-\theta_{1})\tilde{X}_{t_{n+1}%
}^{t_{n}}f_{t_{n+1}}\right)  \Delta t_{n}\right. \label{3.3.1}\\
&  \text{ \  \ }+\left.  \left(  \theta_{2}g_{t_{n}}+(1-\theta_{2})X_{t_{n+1}%
}^{t_{n}}g_{t_{n+1}}\right)  \Delta \langle B\rangle_{n+1}\right]
+R^{n}\nonumber
\end{align}
and
\begin{align}
&  \tilde{X}_{t_{n+1}}^{t_{n}}=\exp \Big(b_{t_{n}}^{11}\Delta B_{n+1}^{1}%
-\frac{1}{2}|b_{t_{n}}^{11}|^{2}\Delta \langle B^{1}\rangle_{n+1}%
\Big)\exp \Big(b_{t_{n}}^{22}\Delta B_{n+1}^{2}-\frac{1}{2}|b_{t_{n}}^{22}%
|^{2}\Delta \langle B^{2}\rangle_{n+1}\Big)\label{3.3.2}\\
&  \exp \Big(a_{t_{n}}^{1}\Delta \tilde{B}_{n+1}^{1}-\frac{1}{2}|a_{t_{n}}%
^{1}|^{2}\Delta \langle \tilde{B}^{1}\rangle_{n+1}\Big)\exp \Big(a_{t_{n}}%
^{2}\Delta \tilde{B}_{n+1}^{2}-\frac{1}{2}|a_{t_{n}}^{2}|^{2}\Delta
\langle \tilde{B}^{2}\rangle_{n+1}\Big)\nonumber \\
&  \exp \Big(b_{t_{n}}^{12}\Delta \tilde{B}_{n+1}^{3}-\frac{1}{2}|b_{t_{n}}%
^{12}|^{2}\Delta \langle \tilde{B}^{3}\rangle_{n+1}\Big)\exp \Big(b_{t_{n}}%
^{21}\Delta \tilde{B}_{n+1}^{4}-\frac{1}{2}|b_{t_{n}}^{21}|^{2}\Delta
\langle \tilde{B}^{4}\rangle_{n+1}\Big)\nonumber \\
&  \exp \Big(-b_{t_{n}}^{12}b_{t_{n}}^{22}\Delta \langle B^{1}\rangle
_{n+1}-b_{t_{n}}^{11}b_{t_{n}}^{21}\Delta \langle B^{2}\rangle_{n+1}%
\Big)\exp \Big(-(a_{t_{n}}^{1}b_{t_{n}}^{11}+a_{t_{n}}^{2}b_{t_{n}}^{22})\Delta
t_{n}\Big),\nonumber
\end{align}
where $\theta_{i}\in \left[  0,1\right]  (i=1,2)$, $g_{t}:=g(t,y)=(g_{11}%
(t,y),g_{22}(t,y))$, $\bigtriangleup \langle B\rangle_{n+1}:=(\bigtriangleup
\langle B^{1}\rangle_{n+1},\bigtriangleup \langle B^{2}\rangle_{n+1})^{\top}$,
and $R^{n}=R_{x}^{n}+R_{y_{1}}^{n}+R_{y_{2}}^{n}$ with
\begin{align*}
R_{y_{1}}^{n}  &  =\mathbb{\hat{E}}_{t_{n}}^{x}\left[  X_{t_{n+1}}^{t_{n}%
}Y_{t_{n+1}}+\int_{t_{n}}^{t_{n+1}}X_{s}^{t_{n}}f\left(  s,Y_{s}\right)
\text{d}s+\int_{t_{n}}^{t_{n+1}}X_{s}^{t_{n}}g_{ii}\left(  s,Y_{s}\right)
d\langle B^{i}\rangle_{s}\right] \\
&  \text{ }-\mathbb{\hat{E}}_{t_{n}}^{x}\left[  X_{t_{n+1}}^{t_{n}}Y_{t_{n+1}%
}+\left(  \theta_{1}f_{t_{n}}+(1-\theta_{1})X_{t_{n+1}}^{t_{n}}f_{t_{n+1}%
}\right)  \Delta t_{n}+\int_{t_{n}}^{t_{n+1}}X_{s}^{t_{n}}g_{ii}\left(
s,Y_{s}\right)  d\langle B^{i}\rangle_{s}\right]  ,
\end{align*}%
\begin{align*}
R_{y_{2}}^{n}  &  =\mathbb{\hat{E}}_{t_{n}}^{x}\left[  X_{t_{n+1}}^{t_{n}%
}Y_{t_{n+1}}+\left(  \theta_{1}f_{t_{n}}+(1-\theta_{1})X_{t_{n+1}}^{t_{n}%
}f_{t_{n+1}}\right)  \Delta t_{n}+\int_{t_{n}}^{t_{n+1}}X_{s}^{t_{n}}%
g_{ii}\left(  s,Y_{s}\right)  d\langle B^{i}\rangle_{s}\right]  \text{ }\\
&  \text{ }-\mathbb{\hat{E}}_{t_{n}}^{x}\left[  X_{t_{n+1}}^{t_{n}}Y_{t_{n+1}%
}+\left(  \theta_{1}f_{t_{n}}+(1-\theta_{1})X_{t_{n+1}}^{t_{n}}f_{t_{n+1}%
}\right)  \Delta t_{n}\right. \\
&  \text{ }+\left.  \left(  \theta_{2}g_{t_{n}}+(1-\theta_{2})X_{t_{n+1}%
}^{t_{n}}g_{t_{n+1}}\right)  \Delta \langle B\rangle_{n+1}\right]  ,
\end{align*}
and%
\begin{align*}
R_{x}^{n}  &  =\mathbb{\hat{E}}_{t_{n}}^{x}\left[  X_{t_{n+1}}^{t_{n}%
}Y_{t_{n+1}}+(\theta_{1}f_{t_{n}}+(1-\theta_{1})X_{t_{n+1}}^{t_{n}}f_{t_{n+1}%
})\Delta t_{n}+(\theta_{2}g_{t_{n}}+(1-\theta_{2})X_{t_{n+1}}^{t_{n}%
}g_{t_{n+1}})\Delta \langle B\rangle_{n+1}\right] \\
&  \text{\ }-\mathbb{\hat{E}}_{t_{n}}^{x}\left[  \tilde{X}_{t_{n+1}}^{t_{n}%
}Y_{t_{n+1}}+(\theta_{1}f_{t_{n}}+(1-\theta_{1})\tilde{X}_{t_{n+1}}^{t_{n}%
}f_{t_{n+1}})\Delta t_{n}+(\theta_{2}g_{t_{n}}+(1-\theta_{2})\tilde
{X}_{t_{n+1}}^{t_{n}}g_{t_{n+1}})\Delta \langle B\rangle_{n+1}\right]  .
\end{align*}
Based on (\ref{3.3.1})-(\ref{3.3.2}), for given $Y^{N}$, we can solve for
$Y^{n}$, $n=N-1,\ldots,0$, as follows%
\[
\left \{
\begin{array}
[c]{l}%
Y^{n}\text{ }=\mathbb{\hat{E}}_{t_{n}}^{x}\left[  X^{n+1}Y^{n+1}+\left(
\theta_{1}f(t_{n},Y^{n})+(1-\theta_{1})X^{n+1}f(t_{n+1},Y^{n+1})\right)
\Delta t_{n}\right. \\
\text{ \  \  \  \ }+\left.  \left(  \theta_{2}g(t_{n},Y^{n})+(1-\theta
_{2})X^{n+1}g(t_{n+1},Y^{n+1})\right)  \Delta \langle B\rangle_{n+1}\right]
,\\
X^{n+1}=\exp \left(  b_{t_{n}}^{11}\Delta B_{n+1}^{1}-\frac{1}{2}|b_{t_{n}%
}^{11}|^{2}\Delta \langle B^{1}\rangle_{n+1}\right)  \exp \left(  b_{t_{n}}%
^{22}\Delta B_{n+1}^{2}-\frac{1}{2}|b_{t_{n}}^{22}|^{2}\Delta \langle
B^{2}\rangle_{n+1}\right) \\
\text{ \  \  \  \ }\exp \left(  a_{t_{n}}^{1}\Delta \tilde{B}_{n+1}^{1}-\frac{1}%
{2}|a_{t_{n}}^{1}|^{2}\Delta \langle \tilde{B}^{1}\rangle_{n+1}\right)
\exp \left(  a_{t_{n}}^{2}\Delta \tilde{B}_{n+1}^{2}-\frac{1}{2}|a_{t_{n}}%
^{2}|^{2}\Delta \langle \tilde{B}^{2}\rangle_{n+1}\right) \\
\text{ \  \  \  \ }\exp \left(  b_{t_{n}}^{12}\Delta \tilde{B}_{n+1}^{3}-\frac
{1}{2}|b_{t_{n}}^{12}|^{2}\Delta \langle \tilde{B}^{3}\rangle_{n+1}\right)
\exp \left(  b_{t_{n}}^{21}\Delta \tilde{B}_{n+1}^{4}-\frac{1}{2}|b_{t_{n}}%
^{21}|^{2}\Delta \langle \tilde{B}^{4}\rangle_{n+1}\right) \\
\text{ \  \  \  \ }\exp \left(  -b_{t_{n}}^{12}b_{t_{n}}^{22}\Delta \langle
B^{1}\rangle_{n+1}-b_{t_{n}}^{11}b_{t_{n}}^{21}\Delta \langle B^{2}%
\rangle_{n+1}\right)  \exp \left(  -(a_{t_{n}}^{1}b_{t_{n}}^{11}+a_{t_{n}}%
^{2}b_{t_{n}}^{22})\Delta t_{n}\right)  ,
\end{array}
\right.
\]
with the deterministic parameters $\theta_{i}\in \left[  0,1\right]  (i=1,2)$.
\end{example}

For a general given sublinear function $G:$ $\mathbb{S}(d)\rightarrow
\mathbb{R}$ such that%
\begin{equation}
G\left(  A\right)  =\frac{1}{2}\sup_{Q\in \Sigma}tr[AQ], \label{GA}%
\end{equation}
with the bounded and closed set $\Sigma \subset$ $\mathbb{S}_{+}(d)$, and
$B_{t}=(B_{t}^{1},\ldots,B_{t}^{d})^{\top}$ be the $d$-dimensional
$G$-Brownian motion on the $G$-expectation space $(\Omega_{T},L_{G}^{1}(\Omega_{T}),\mathbb{\hat{E}})$ 
with $\Omega_{T}=C_{0}([0,T];\mathbb{R}^{d})$. 
In order to get the explicit expression of $Y_{\cdot}$, we can
construct a suitable auxiliary space $(\tilde{\Omega}_{T},L_{\tilde{G}}%
^{1}(\tilde{\Omega}_{T}),\mathbb{\hat{E}})$ with $\tilde{\Omega}_{T}%
=C_{0}([0,T],\mathbb{R}^{d+d^{2}})$ and define the canonical process
$(B_{t}^{1},\ldots,B_{t}^{d},\tilde{B}_{t}^{1},\ldots,\tilde{B}_{t}^{d^{2}%
})^{\top}$ on this extended space. By introducing the following $\tilde{G}%
$-SDE on $[0,T]$,
\[
X_{t}=1+\sum_{i=1}^{d}\int_{0}^{t}C_{s}^{i}X_{s}dB_{s}^{i}+\sum_{j=1}^{d^{2}%
}\int_{0}^{t}C_{s}^{d+j}X_{s}d\tilde{B}_{s}^{j}%
\]
with some selected coefficients $C_{\cdot}^{i}$, $i=1,\ldots,d+d^{2}$,
$G$-BSDE (\ref{G-BSDE}) can be solved explicitly by%
\[
Y_{t}=\mathbb{\hat{E}}_{t}\left[  X_{T}^{t}\xi+\int_{t}^{T}X_{s}^{t}f\left(
s,Y_{s}\right)  \text{d}s+\int_{t}^{T}X_{s}^{t}g_{ij}\left(  s,Y_{s}\right)
d\langle B^{i},B^{j}\rangle_{s}\right]  \text{,}%
\]
with $X_{s}^{t}=X_{s}/X_{t}$. Similar to the previous processes, we get the
reference equation%
\begin{align}
Y_{t_{n}}  &  =\mathbb{\hat{E}}_{t_{n}}^{x}\left[  \tilde{X}_{t_{n+1}}^{t_{n}%
}Y_{t_{n+1}}+\left(  \theta_{1}f_{t_{n}}+(1-\theta_{1})\tilde{X}_{t_{n+1}%
}^{t_{n}}f_{t_{n+1}}\right)  \Delta t_{n}\right. \label{multi-reference}\\
&  \text{ \  \ }+\left.  \left(  \theta_{2}g_{t_{n}}+(1-\theta_{2})X_{t_{n+1}%
}^{t_{n}}g_{t_{n+1}},\Delta \langle B\rangle_{n+1}\right)  \right]
+R^{n},\nonumber
\end{align}
where $\bigtriangleup \langle B\rangle_{n+1}:=(\bigtriangleup \langle
B^{i},B^{j}\rangle_{n+1})_{i,j=1}^{d}$,
\begin{equation}
\tilde{X}_{t_{n+1}}^{t_{n}}=1+\sum_{i=1}^{d}\int_{t_{n}}^{t_{n+1}}C_{t_{n}%
}^{i}\tilde{X}_{s}^{t_{n}}dB_{s}^{i}+\sum_{j=1}^{d^{2}}\int_{t_{n}}^{t_{n+1}%
}C_{t_{n}}^{d+j}\tilde{X}_{s}^{t_{n}}d\tilde{B}_{s}^{j}, \label{G-SDE-multi}%
\end{equation}
and $R^{n}$ is the truncation error in the approximation of the forward and
backward equations.

Based on (\ref{multi-reference}), we propose the following $\theta$-scheme for
solving $Y_{t}$ of $G$-BSDE (\ref{G-BSDE}).
\begin{sch}
\label{Scheme multi-dim} Given random variables $Y^{N}$, solve random
variables $Y^{n}$, $n=N-1,\ldots,0$ from%
\begin{equation}
\left \{
\begin{array}
[c]{l}%
Y^{n}\text{ }=\mathbb{\hat{E}}_{t_{n}}^{x}\left[  X^{n+1}Y^{n+1}+\left(
\theta_{1}f(t_{n},Y^{n})+(1-\theta_{1})X^{n+1}f(t_{n+1},Y^{n+1})\right)
\Delta t_{n}\right. \\
\text{ \  \  \  \  \ }+\left.  \left(  \theta_{2}g(t_{n},Y^{n})+(1-\theta
_{2})X^{n+1}g(t_{n+1},Y^{n+1}),\Delta \langle B\rangle_{n+1}\right)  \right]
,\\
X^{n+1}\text{ is the solution of }\tilde{G}\text{-SDE }%
\eqref{G-SDE-multi}\text{ at }t=t_{n+1},
\end{array}
\right.  \label{multi-scheme}%
\end{equation}
with the deterministic parameters $\theta_{i}\in \left[  0,1\right]  (i=1,2)$.
\end{sch}

\section{Convergence analysis of the discrete-time scheme}

In this section, we focus on the convergence analysis of the discrete
approximation scheme. For the reader's convenience, the proofs are done with
$d=1$, and the results still hold for the case $d>1$.

\subsection{A stability theorem}

We first present a useful theorem in our convergence analysis. In what
follows, $C$ represents a generic constant, which may be different from line
to line.

\begin{lemma}
\label{lemma1}For any $t\in \lbrack t_{n},t_{n+1}]$, $n=0,1,\ldots,N-1$, we have
\begin{description}
\item[(i)] $\mathbb{\hat{E}}_{t_{n}}^{x}[\tilde{X}_{t}^{t_{n}}]=1$ and
$\mathbb{\hat{E}}_{t_{n}}^{x}[(\tilde{X}_{t}^{t_{n}})^{2}]\leq1+C\Delta t$;

\item[(ii)] $\mathbb{\hat{E}}_{t_{n}}^{x}[X_{t}^{t_{n}}]=1$ and $\mathbb{\hat
{E}}_{t_{n}}^{x}[(X_{t}^{t_{n}})^{2}]\leq1+C\Delta t$,
\end{description}
where $C>0$ is a constant independent of $n$, $X_{\cdot}^{t_{n}}$, and
$\tilde{X}_{\cdot}^{t_{n}}$.

\begin{proof}
(i) For any $t\in \lbrack t_{n},t_{n+1}]$, from (\ref{3.2.4}) we know%
\begin{align*}
\tilde{X}_{t}^{t_{n}}  &  =\exp \left(  -\int_{t_{n}}^{t}a_{t_{n}}b_{t_{n}%
}ds\right)  \exp \left(  \int_{t_{n}}^{t}b_{t_{n}}dB_{s}-\frac{1}{2}\int
_{t_{n}}^{t}b_{t_{n}}^{2}d\langle B\rangle_{s}\right) \\
&  \text{ \  \ }\exp \left(  \int_{t_{n}}^{t}a_{t_{n}}d\tilde{B}_{s}-\frac{1}%
{2}\int_{t_{n}}^{t}a_{t_{n}}^{2}d\langle \tilde{B}\rangle_{s}\right)  .
\end{align*}
Notice that $\langle B,\tilde{B}\rangle_{t}=t$, from $G$-It\^{o}'s formula, we
get
\begin{equation}
\tilde{X}_{t}^{t_{n}}=1+\int_{t_{n}}^{t}a_{t_{n}}\tilde{X}_{s}^{t_{n}}%
d\tilde{B}_{s}+\int_{t_{n}}^{t}b_{t_{n}}\tilde{X}_{s}^{t_{n}}dB_{s}.
\label{3.1}%
\end{equation}
From Proposition \ref{PropE} (ii), we have $\mathbb{\hat{E}}_{t_{n}}%
^{x}[\tilde{X}_{t}^{t_{n}}]=1$. Since $\langle B\rangle_{t}-\langle
B\rangle_{t_{n}}\leq(t-t_{n})$ and $\langle \tilde{B}\rangle_{t}-\langle
\tilde{B}\rangle_{t_{n}}\leq \frac{1}{\sigma^{2}}(t-t_{n})$, for any
$t\in \lbrack t_{n},t_{n+1}]$, it follows that%
\begin{align*}
(\tilde{X}_{t}^{t_{n}})^{2}  &  =\exp \left(  \int_{t_{n}}^{t}a_{t_{n}}%
^{2}d\langle \tilde{B}\rangle_{s}+\int_{t_{n}}^{t}b_{t_{n}}^{2}d\langle
B\rangle_{s}+2\int_{t_{n}}^{t}a_{t_{n}}b_{t_{n}}ds\right)  \exp \left(
-4\int_{t_{n}}^{t}a_{t_{n}}b_{t_{n}}ds\right) \\
&  \text{ \  \ }\exp \left(  2\int_{t_{n}}^{t}a_{t_{n}}d\tilde{B}_{s}%
-2\int_{t_{n}}^{t}a_{t_{n}}^{2}d\langle \tilde{B}\rangle_{s}\right)
\exp \left(  2\int_{t_{n}}^{t}b_{t_{n}}dB_{s}-2\int_{t_{n}}^{t}b_{t_{n}}%
^{2}d\langle B\rangle_{s}\right) \\
&  \leq \exp \left(  (a_{t_{n}}+b_{t_{n}})^{2}\Delta t+(\frac{1}{\sigma^{2}%
}-1)a_{t_{n}}^{2}\Delta t\right)  R_{t}^{t_{n}},
\end{align*}
where
\[
R_{t}^{t_{n}}:=\exp \left(  -4\int_{t_{n}}^{t}a_{t_{n}}b_{t_{n}}ds\right)
\  \exp \left(  2\int_{t_{n}}^{t}a_{t_{n}}d\tilde{B}_{s}-2\int_{t_{n}}%
^{t}a_{t_{n}}^{2}d\langle \tilde{B}\rangle_{s}\right)  \exp \left(  2\int
_{t_{n}}^{t}b_{t_{n}}dB_{s}-2\int_{t_{n}}^{t}b_{t_{n}}^{2}d\langle
B\rangle_{s}\right)  .
\]
Applying $G$-It\^{o}'s formula to $R_{t}^{t_{n}}$, we can deduce
$\mathbb{\hat{E}}_{t_{n}}^{x}[R_{t}^{t_{n}}]=1$ for $t\in \lbrack t_{n}%
,t_{n+1}]$. This implies that
\[
\mathbb{\hat{E}}_{t_{n}}^{x}[(\tilde{X}_{t}^{t_{n}})^{2}]\leq \exp \left(
\left[  (a_{t_{n}}+b_{t_{n}})^{2}+(\frac{1}{\sigma^{2}}-1)a_{t_{n}}%
^{2}\right]  \Delta t\right)  .
\]
Noting that $e^{x}\leq1+ex$, for any $0\leq x\leq1$, then for bounded
processes $a_{s}$, $b_{s}$, and sufficiently small $\Delta t$, we draw the
conclusion (i). Similarly, we can prove (ii).
\end{proof}
\end{lemma}

\begin{theorem}
\label{theorem1}Assume that $f$ and $g$ satisfy (H1)-(H2) with Lipschitz
constant $L.$ Then for sufficiently small time step $\Delta t,$ we have%
\begin{equation}
\mathbb{\hat{E}[}|Y_{t_{n}}-Y^{n}|^{2}]\  \  \leq C\left(  \mathbb{\hat{E}%
[}|Y_{T}-Y^{N}|^{2}]+\sum \limits_{i=0}^{N-1}\frac{\mathbb{\hat{E}[}|R_{x}%
^{i}|^{2}]+\mathbb{\hat{E}[}|R_{y}^{i}|^{2}]}{\Delta t}\right)
,\  \label{theorem4.1}%
\end{equation}
for $n=N-1,\ldots,0$, where $C>0$ is a constant just depending on $T$, $L$ and
the upper bounds of $a$ and $b$, and $R_{x}^{n}$ and $R_{y}^{n}$ are error
terms defined in (\ref{RX}) and (\ref{RY1})-(\ref{RY2}), respectively.

\begin{proof}
Let $e_{y}^{n}=Y_{t_{n}}-Y^{n}$, $e_{f}^{n}=f(t_{n},Y_{t_{n}})-f(t_{n},Y^{n})$
and $e_{g}^{n}=g(t_{n},Y_{t_{n}})-g(t_{n},Y^{n})$, for $n=0,1,\ldots,N$.
Subtract $\left(  \ref{scheme1}\right)  $ from $\left(  \ref{reference}%
\right)  $, by Proposition \ref{PropE} $(iii)$, we get
\begin{equation}%
\begin{array}
[c]{l}%
\left \vert e_{y}^{n}\right \vert \leq \mathbb{\hat{E}}_{t_{n}}^{x}\left[
\tilde{X}_{t_{n+1}}^{t_{n}}|e_{y}^{n+1}|+\left(  \theta_{1}|e_{f}%
^{n}|+(1-\theta_{1})\tilde{X}_{t_{n+1}}^{t_{n}}|e_{f}^{n+1}|\right)  \Delta
t_{n}\right. \\
\text{ \  \  \  \  \  \ }+\left.  \left(  \theta_{2}|e_{g}^{n}|+(1-\theta
_{2})\tilde{X}_{t_{n+1}}^{t_{n}}|e_{g}^{n+1}|\right)  \Delta \langle
B\rangle_{n,t_{n+1}}\right]  +|R_{x}^{n}|+|R_{y}^{n}|.
\end{array}
\label{4.1}%
\end{equation}
From $\langle B\rangle_{t_{n+1}}-\langle B\rangle_{t_{n}}\leq \Delta t_{n}$, we
get
\begin{equation}
\left \vert e_{y}^{n}\right \vert \leq(1+2L\Delta t_{n})\mathbb{\hat{E}}_{t_{n}%
}^{x}[\tilde{X}_{t_{n+1}}^{t_{n}}|e_{y}^{n+1}|]+2L\Delta t_{n}\left \vert
e_{y}^{n}\right \vert +|R_{x}^{n}|+|R_{y}^{n}|, \label{4.2}%
\end{equation}
where $L>0$ is the Lipschitz constant. This implies that
\begin{align}
\left \vert e_{y}^{n}\right \vert  &  \leq \frac{1+2L\Delta t_{n}}{1-2L\Delta
t_{n}}\mathbb{\hat{E}}_{t_{n}}^{x}[\tilde{X}_{t_{n+1}}^{t_{n}}|e_{y}%
^{n+1}|]+\frac{1}{1-2L\Delta t_{n}}\left(  |R_{x}^{n}|+|R_{y}^{n}|\right)
\label{4.3}\\
&  \leq(1+C\Delta t)\mathbb{\hat{E}}_{t_{n}}^{x}[\tilde{X}_{t_{n+1}}^{t_{n}%
}|e_{y}^{n+1}|]+C\left(  |R_{x}^{n}|+|R_{y}^{n}|\right)  .\nonumber
\end{align}
By means of the inequalities $(a+b)^{2}\leq(1+\Delta t)a^{2}+(1+\frac
{1}{\Delta t})b^{2}$, H\"{o}lder's inequality, and Lemma \ref{lemma1}, we can
derive that
\begin{equation}
|e_{y}^{n}|^{2}\leq(1+C\Delta t)\mathbb{\hat{E}}_{t_{n}}^{x}[|e_{y}^{n+1}%
|^{2}]+C(1+\frac{1}{\Delta t})\left(  |R_{x}^{n}|^{2}+|R_{y}^{n}|^{2}\right)
. \label{4.4}%
\end{equation}
Taking $\tilde{G}$-expectation on both sides of $\left(  \ref{4.4}\right)  $
and substituting $e_{y}^{i}$, $i=n+1,\ldots,N-1$ recursively, we consequently
deduce that%
\begin{equation}%
\begin{array}
[c]{ll}%
\mathbb{\hat{E}[}|e_{y}^{n}|^{2}] & \leq(1+C\Delta t)^{N-n}\mathbb{\hat{E}%
}[|e_{y}^{N}|^{2}]+\frac{C}{\Delta t}\sum \limits_{i=n}^{N-1}(1+\Delta
t)^{i-n}\left(  \mathbb{\hat{E}[}|R_{x}^{i}|^{2}]+\mathbb{\hat{E}[}|R_{y}%
^{i}|^{2}]\right) \\
& \leq e^{C(N-n)\Delta t}\mathbb{\hat{E}}[|e_{y}^{N}|^{2}]+Ce^{(N-n)\Delta
t}\sum \limits_{i=n}^{N-1}\frac{1}{\Delta t}\left(  \mathbb{\hat{E}[}|R_{x}%
^{i}|^{2}]+\mathbb{\hat{E}[}|R_{y}^{i}|^{2}]\right) \\
& \leq C\mathbb{\hat{E}}[|e_{y}^{N}|^{2}]+C\sum \limits_{i=n}^{N-1}\frac
{1}{\Delta t}\left(  \mathbb{\hat{E}[}|R_{x}^{i}|^{2}]+\mathbb{\hat{E}[}%
|R_{y}^{i}|^{2}]\right)  ,
\end{array}
\label{4.5}%
\end{equation}
which concludes the proof.
\end{proof}
\end{theorem}

\subsection{Error estimates}

Next, we will give the estimates of truncation errors $R_{x}^{n}$ and
$R_{y}^{n}$ defined in (\ref{RX}) and (\ref{RY1})-(\ref{RY2}). Then we come to
the conclusion based on Theorem \ref{theorem1}.

\begin{lemma}
\label{lemma2}Let $R_{x}^{n}$ be the local truncation error derived from
(\ref{RX}). Then for sufficiently small time step $\Delta t$ and $\theta
_{i}\left(  i=1,2\right)  \in \lbrack0,1]$, we have the following conclusions:

\begin{description}
\item[(i)] If $a_{s}=a(s,B_{s})$, $b_{s}=b(s,B_{s})\in C^{\frac{1}{2},1}$ are
bounded stochastic processes on $[0,T]$ and $f,g,\varphi \in C^{\frac{1}{2},1}%
$, we have%
\begin{equation}
\left \vert R_{x}^{n}\right \vert \leq C(\Delta t)^{\frac{3}{2}}.
\label{truncation RX1}%
\end{equation}

\item[(ii)] In particular, if $a_{s}=a(s)$, $b_{s}=b(s)$ are Lipschitz
continuous functions on $[0,T]$ and $f,g,u\in C_{b}^{1,2}$, we have%
\begin{equation}
\left \vert R_{x}^{n}\right \vert \leq C(\Delta t)^{2}. \label{truncation RX2}%
\end{equation}
\end{description}
Here $C>0$ is a constant just depending on $T$, $L$, and the upper bounds of
$a$ and $b$.

\begin{proof}
(i) From the definition of $R_{x}^{n}$ and Proposition \ref{PropE} $(iii)$, we
can deduce%
\begin{align}
|R_{x}^{n}|  &  \leq \mathbb{\hat{E}}_{t_{n}}^{x}\big[(X_{t_{n+1}}^{t_{n}}
-\tilde{X}_{t_{n+1}}^{t_{n}})(Y_{t_{n+1}}+f_{t_{n+1}}\Delta t_{n}+g_{t_{n+1}
}\Delta \langle B\rangle_{n,t_{n+1}})\big]\nonumber \\
&  \text{ \  \ }\vee \mathbb{\hat{E}}_{t_{n}}^{x}\big[(\tilde{X}_{t_{n+1}}^{t_{n}
}-X_{t_{n+1}}^{t_{n}})(Y_{t_{n+1}}+f_{t_{n+1}}\Delta t_{n}+g_{t_{n+1}}
\Delta \langle B\rangle_{n,t_{n+1}})\big]\nonumber \\
&  \leq \left \{  \mathbb{\hat{E}}_{t_{n}}^{x}\big[(X_{t_{n+1}}^{t_{n}}-\tilde
{X}_{t_{n+1}}^{t_{n}})Y_{t_{n+1}}\big]+\mathbb{\hat{E}}_{t_{n}}^{x}\big[(X_{t_{n+1}
}^{t_{n}}-\tilde{X}_{t_{n+1}}^{t_{n}})f_{t_{n+1}}\Delta t_{n}\big]\right.
\label{5.10}\\
&  \text{ \  \ }+\left.  \mathbb{\hat{E}}_{t_{n}}^{x}\big[(X_{t_{n+1}}^{t_{n}
}-\tilde{X}_{t_{n+1}}^{t_{n}})g_{t_{n+1}}\Delta \langle B\rangle_{n,t_{n+1}
}\big]\right \} \nonumber \\
&  \text{ \  \ }\vee \left \{  \mathbb{\hat{E}}_{t_{n}}^{x}\big[(\tilde{X}_{t_{n+1}
}^{t_{n}}-X_{t_{n+1}}^{t_{n}})Y_{t_{n+1}}\big]+\mathbb{\hat{E}}_{t_{n}}
^{x}\big[(\tilde{X}_{t_{n+1}}^{t_{n}}-X_{t_{n+1}}^{t_{n}})f_{t_{n+1}}\Delta
t_{n}\big]\right. \nonumber \\
&  \text{ \  \ }+\left.  \mathbb{\hat{E}}_{t_{n}}^{x}\big[(\tilde{X}_{t_{n+1}
}^{t_{n}}-X_{t_{n+1}}^{t_{n}})g_{t_{n+1}}\Delta \langle B\rangle_{n,t_{n+1}
}\big]\right \}  .\nonumber
\end{align}
We just give the estimate of $\mathbb{\hat{E}}_{t_{n}}^{x}[(X_{t_{n+1}}
^{t_{n}}-\tilde{X}_{t_{n+1}}^{t_{n}})Y_{t_{n+1}}]$, and $\mathbb{\hat{E}
}_{t_{n}}^{x}[(\tilde{X}_{t_{n+1}}^{t_{n}}-X_{t_{n+1}}^{t_{n}})Y_{t_{n+1}}]$
can be similarly obtained. According to $\left(  \ref{3.2.5}\right)  $,
$\left(  \ref{3.1}\right)  $ and Proposition \ref{PropE} (ii), it is worth
noting that%
\begin{align*}
0  &  =\mathbb{\hat{E}}_{t_{n}}^{x}\big[\pm Y_{t_{n}}(X_{t_{n+1}}^{t_{n}}
-\tilde{X}_{t_{n+1}}^{t_{n}})\big]\\
&  =Y_{t_{n}}^{+}\mathbb{\hat{E}}_{t_{n}}^{x}\big[\pm(X_{t_{n+1}}^{t_{n}}
-\tilde{X}_{t_{n+1}}^{t_{n}})\big]+Y_{t_{n}}^{-}\mathbb{\hat{E}}_{t_{n}}^{x}
\big[\mp(X_{t_{n+1}}^{t_{n}}-\tilde{X}_{t_{n+1}}^{t_{n}})\big]\\
&  =Y_{t_{n}}^{+}\mathbb{\hat{E}}_{t_{n}}^{x}\bigg[\pm \int_{t_{n}}^{t_{n+1}
}(a_{t}X_{t}^{t_{n}}-a_{t_{n}}\tilde{X}_{t}^{t_{n}})d\tilde{B}_{t}\pm
\int_{t_{n}}^{t_{n+1}}(b_{t}X_{t}^{t_{n}}-b_{t_{n}}\tilde{X}_{t}^{t_{n}
})dB_{t}\bigg] \\
&  \text{ \  \ }+Y_{t_{n}}^{-}\mathbb{\hat{E}}_{t_{n}}^{x}\bigg[  \mp
\int_{t_{n}}^{t_{n+1}}(a_{t_{n}}\tilde{X}_{t}^{t_{n}}-a_{t}X_{t}^{t_{n}
})d\tilde{B}_{t}\mp \int_{t_{n}}^{t_{n+1}}(b_{t_{n}}\tilde{X}_{t}^{t_{n}}
-b_{t}X_{t}^{t_{n}})dB_{t}\bigg]
\end{align*}
This yields that
\begin{equation}
\mathbb{\hat{E}}_{t_{n}}^{x}\big[(X_{t_{n+1}}^{t_{n}}-\tilde{X}_{t_{n+1}}^{t_{n}%
})Y_{t_{n+1}}\big]=\mathbb{\hat{E}}_{t_{n}}^{x}\big[(X_{t_{n+1}}^{t_{n}}-\tilde
{X}_{t_{n+1}}^{t_{n}})(Y_{t_{n+1}}-Y_{t_{n}})\big]. \label{5.0}%
\end{equation}
Under the conditions $f,g,\varphi \in C^{\frac{1}{2},1}$, from Theorem
\ref{Nonlinear Feynman-Kac}, it is easy to check that $Y_{t}=u(t,B_{t})\in
C^{\frac{1}{2},1}$. Then, it follows from H\"{o}lder's inequality
\begin{equation}
\mathbb{\hat{E}}_{t_{n}}^{x}[(X_{t_{n+1}}^{t_{n}}-\tilde{X}_{t_{n+1}}^{t_{n}%
})(Y_{t_{n+1}}-Y_{t_{n}})]\leq C\left(  \mathbb{\hat{E}}_{t_{n}}%
^{x}[|X_{t_{n+1}}^{t_{n}}-\tilde{X}_{t_{n+1}}^{t_{n}}|^{2}]\right)  ^{\frac
{1}{2}}\sqrt{\Delta t_{n}}. \label{5.9}%
\end{equation}
Note that%
\begin{equation}
\mathbb{\hat{E}}_{t_{n}}^{x}\left[  \left \vert \int_{t_{n}}^{t_{n+1}}\eta
_{t}d\langle B\rangle_{t}\right \vert \right]  \leq \mathbb{\hat{E}}_{t_{n}}%
^{x}\left[  \int_{t_{n}}^{t_{n+1}}\left \vert \eta_{t}\right \vert dt\right]
,\text{ }\eta \in M_{G}^{1}(0,T), \label{<B><=B}%
\end{equation}
and%
\begin{equation}
\text{ \  \ }\mathbb{\hat{E}}_{t_{n}}^{x}\left[  \left(  \int_{t_{n}}^{t_{n+1}%
}\eta_{t}dB_{t}\right)  ^{2}\right]  \leq \mathbb{\hat{E}}_{t_{n}}^{x}\left[
\int_{t_{n}}^{t_{n+1}}|\eta_{t}|^{2}d\langle B\rangle_{t}\right]  ,\text{
}\eta \in M_{G}^{2}(0,T). \label{isometry}%
\end{equation}
%Applying The Euler Approximation To (\Ref{<B><=B}) And (\Ref{Isometry}) (Cf.
%Theorem 10.2.2 In \Cite{Kp1992}), We Can Get
Combining (\ref{3.2.5}), (\ref{3.1}), (\ref{<B><=B}) and (\ref{isometry}), we
have
\begin{align}
&  \mathbb{\hat{E}}_{t_{n}}^{x}[|X_{t}^{t_{n}}-\tilde{X}_{t}^{t_{n}}%
|^{2}]\nonumber \\
&  \leq2\mathbb{\hat{E}}_{t_{n}}^{x}\left[  \int_{t_{n}}^{t}(a_{s}X_{s}%
^{t_{n}}-a_{t_{n}}\tilde{X}_{s}^{t_{n}})^{2}d\langle \tilde{B}\rangle
_{s}\right]  +2\mathbb{\hat{E}}_{t_{n}}^{x}\left[  \int_{t_{n}}^{t}(b_{s}%
X_{s}^{t_{n}}-b_{t_{n}}\tilde{X}_{s}^{t_{n}})^{2}d\langle B\rangle_{s}\right]
\nonumber \\
&  \leq2\int_{t_{n}}^{t}\mathbb{\hat{E}}_{t_{n}}^{x}\left[  \left(
a_{s}(X_{s}^{t_{n}}-\tilde{X}_{s}^{t_{n}})+(a_{s}-a_{t_{n}})\tilde{X}%
_{s}^{t_{n}}\right)  ^{2}\right]  ds\\
&  \text{ \  \ }+2\int_{t_{n}}^{t}\mathbb{\hat{E}}_{t_{n}}^{x}\left[  \left(
b_{s}(X_{s}^{t_{n}}-\tilde{X}_{s}^{t_{n}})+(b_{s}-b_{t_{n}})\tilde{X}%
_{s}^{t_{n}}\right)  ^{2}\right]  ds\nonumber \\
&  \leq C\int_{t_{n}}^{t}\mathbb{\hat{E}}_{t_{n}}^{x}[(X_{s}^{t_{n}}-\tilde
{X}_{s}^{t_{n}})^{2}]ds+C\int_{t_{n}}^{t}\mathbb{\hat{E}}_{t_{n}}^{x}%
[(a_{s}-a_{t_{n}})^{2}(\tilde{X}_{s}^{t_{n}})^{2}]ds\nonumber \\
&  \text{ \  \ }+C\int_{t_{n}}^{t}\mathbb{\hat{E}}_{t_{n}}^{x}[(b_{s}-b_{t_{n}%
})^{2}(\tilde{X}_{s}^{t_{n}})^{2}]ds.\nonumber
\end{align}
Similar to Lemma \ref{lemma1}, we can prove $\mathbb{\hat{E}}_{t_{n}}%
^{x}[(\tilde{X}_{t}^{t_{n}})^{4}]\leq C<\infty$. Using H\"{o}lder's inequality
again, we obtain%
\begin{align}
&  \mathbb{\hat{E}}_{t_{n}}^{x}[|X_{t}^{t_{n}}-\tilde{X}_{t}^{t_{n}}%
|^{2}]\nonumber \\
&  \leq C\int_{t_{n}}^{t}\mathbb{\hat{E}}_{t_{n}}^{x}[(X_{s}^{t_{n}}-\tilde
{X}_{s}^{t_{n}})^{2}]ds+C\int_{t_{n}}^{t}\left(  \mathbb{\hat{E}}_{t_{n}}%
^{x}[(a_{s}-a_{t_{n}})^{4}]\right)  ^{\frac{1}{2}}\left(  \mathbb{\hat{E}%
}_{t_{n}}^{x}[(\tilde{X}_{s}^{t_{n}})^{4}]\right)  ^{\frac{1}{2}}ds\nonumber \\
&  \text{ \  \ }+C\int_{t_{n}}^{t}\left(  \mathbb{\hat{E}}_{t_{n}}^{x}%
[(b_{s}-b_{t_{n}})^{4}]\right)  ^{\frac{1}{2}}\left(  \mathbb{\hat{E}}_{t_{n}%
}^{x}[(\tilde{X}_{s}^{t_{n}})^{4}]\right)  ^{\frac{1}{2}}ds\label{5.2}\\
&  \leq C\int_{t_{n}}^{t}\mathbb{\hat{E}}_{t_{n}}^{x}[(X_{s}^{t_{n}}-\tilde
{X}_{s}^{t_{n}})^{2}]ds+C\int_{t_{n}}^{t}\left(  \mathbb{\hat{E}}_{t_{n}}%
^{x}[(a_{s}-a_{t_{n}})^{4}]\right)  ^{\frac{1}{2}}ds\nonumber \\
&  \text{ \  \ }+C\int_{t_{n}}^{t}\left(  \mathbb{\hat{E}}_{t_{n}}^{x}%
[(b_{s}-b_{t_{n}})^{4}]\right)  ^{\frac{1}{2}}ds.\nonumber
\end{align}
Note that $\mathbb{\hat{E}}[|B_{s}-B_{t_{n}}|^{2k}]=3\bar{\sigma}^{4}%
(s-t_{n})^{k}$, for any $k\in \mathbb{N}^{+}$, If $a_{s}$ $=a(s,B_{s})$ and
$b_{s}=b(s,B_{s})$ are bounded stochastic processes in $C^{\frac{1}{2},1}$, we
then get
\begin{equation}%
\begin{array}
[c]{rr}%
\mathbb{\hat{E}}_{t_{n}}^{x}[(a_{s}-a_{t_{n}})^{4}]\leq C(s-t_{n})^{2}\text{,}
& \mathbb{\hat{E}}_{t_{n}}^{x}[(b_{s}-b_{t_{n}})^{4}]\leq C(s-t_{n}%
)^{2}\text{.}%
\end{array}
\label{5.4}%
\end{equation}
Combining (\ref{5.2}) and (\ref{5.4}), from Gronwall's inequality, we have
\begin{equation}
\mathbb{\hat{E}}_{t_{n}}^{x}[|X_{t}^{t_{n}}-\tilde{X}_{t}^{t_{n}}|^{2}]\leq
C(t-t_{n})^{2}. \label{5.6}%
\end{equation}
Then it follows from (\ref{5.0}), (\ref{5.9}), and (\ref{5.6}) that
\[
\mathbb{\hat{E}}_{t_{n}}^{x}[\pm(X_{t_{n+1}}^{t_{n}}-\tilde{X}_{t_{n+1}%
}^{t_{n}})Y_{t_{n+1}}]\leq C(\Delta t_{n})^{\frac{3}{2}}.
\]
In the same way above, we can also obtain $\mathbb{\hat{E}}_{t_{n}}^{x}%
[\pm(X_{t_{n+1}}^{t_{n}}-\tilde{X}_{t_{n+1}}^{t_{n}})f_{t_{n+1}}\Delta
t_{n}]\leq C(\Delta t_{n})^{\frac{5}{2}}$ and $\mathbb{\hat{E}}_{t_{n}}%
^{x}[\pm(X_{t_{n+1}}^{t_{n}}-\tilde{X}_{t_{n+1}}^{t_{n}})g_{t_{n+1}}%
\Delta \langle B\rangle_{n,t_{n+1}}]\leq C(\Delta t_{n})^{\frac{5}{2}}$. Thus,
(\ref{truncation RX1}) follows from (\ref{5.10}).

(ii) In particular, if $f,$ $g,$ $u\in C_{b}^{1,2}$, applying $G$-It\^{o}'s
formula to $(X_{t}^{t_{n}}-\tilde{X}_{t}^{t_{n}})(Y_{t}-Y_{t_{n}})$ on
$[t_{n,}$ $t_{n+1}]$ and taking the conditional $G$-expectation $\mathbb{\hat
{E}}_{t_{n}}^{x}[\cdot]$, we have%
\begin{align}
&  \mathbb{\hat{E}}_{t_{n}}^{x}[(X_{t_{n+1}}^{t_{n}}-\tilde{X}_{t_{n+1}%
}^{t_{n}})(Y_{t_{n+1}}-Y_{t_{n}})]\\
&  =\mathbb{\hat{E}}_{t_{n}}^{x}\left[  \int_{t_{n}}^{t_{n+1}}(X_{t}^{t_{n}%
}-\tilde{X}_{t}^{t_{n}})\left(  \partial_{t}u(t,B_{t})+a_{t}\partial
_{x}u(t,B_{t})\right)  dt\right. \nonumber \\
&  \text{ \  \ }+\int_{t_{n}}^{t_{n+1}}(X_{t}^{t_{n}}-\tilde{X}_{t}^{t_{n}%
})(\frac{1}{2}\partial_{xx}^{2}u(t,B_{t})+b_{t}\partial_{x}u(t,B_{t}))d\langle
B\rangle_{t}\nonumber \\
&  \text{ \  \ }+\left.  \int_{t_{n}}^{t_{n+1}}(a_{t}-a_{t_{n}})\partial
_{x}u(t,B_{t})\tilde{X}_{t}^{t_{n}}dt+\int_{t_{n}}^{t_{n+1}}(b_{t}-b_{t_{n}%
})\partial_{x}u(t,B_{t})\tilde{X}_{t}^{t_{n}}d\langle B\rangle_{t}\right]
.\nonumber
\end{align}
By means of (\ref{<B><=B}) and H\"{o}lder's inequality, we can deduce
\begin{align}
&  \mathbb{\hat{E}}_{t_{n}}^{x}[(X_{t_{n+1}}^{t_{n}}-\tilde{X}_{t_{n+1}%
}^{t_{n}})(Y_{t_{n+1}}-Y_{t_{n}})]\label{5.7}\\
&  \leq C\int_{t_{n}}^{t_{n+1}}\left(  \mathbb{\hat{E}}_{t_{n}}^{x}%
[|X_{t}^{t_{n}}-\tilde{X}_{t}^{t_{n}}|^{2}]\right)  ^{\frac{1}{2}%
}dt\nonumber \\
&  \text{ \  \ }+C\int_{t_{n}}^{t_{n+1}}\left(  \mathbb{\hat{E}}_{t_{n}}%
^{x}[|a_{t}-a_{t_{n}}|^{2}]\right)  ^{\frac{1}{2}}\left(  \mathbb{\hat{E}%
}_{t_{n}}^{x}[(\tilde{X}_{t}^{t_{n}})^{2}]\right)  ^{\frac{1}{2}}dt\nonumber \\
&  \text{ \  \ }+C\int_{t_{n}}^{t_{n+1}}\left(  \mathbb{\hat{E}}_{t_{n}}%
^{x}[|b_{t}-b_{t_{n}}|^{2}]\right)  ^{\frac{1}{2}}\left(  \mathbb{\hat{E}%
}_{t_{n}}^{x}[(\tilde{X}_{t}^{t_{n}})^{2}]\right)  ^{\frac{1}{2}}dt.\nonumber
\end{align}
Since $a_{s}=a(s)$, $b_{s}=b(s)$ are Lipschitz continuous functions on
$[0,T]$, from (\ref{5.2}) and (\ref{5.7}), we can easily obtain%
\[
\mathbb{\hat{E}}_{t_{n}}^{x}[\pm(X_{t_{n+1}}^{t_{n}}-\tilde{X}_{t_{n+1}%
}^{t_{n}})Y_{t_{n+1}}]\leq C(\Delta t_{n})^{2}.
\]
Similar to (\ref{truncation RX1}), (\ref{truncation RX2}) holds.
\end{proof}
\end{lemma}

\begin{remark}
Note that Lemma \ref{lemma2} (ii) requires smoothness condition on $u$. In the
following, we give some assumptions that implies the needed smoothness for
$u$. For given constants $\alpha,\beta \in(0,1]$, denote
\[
\left \Vert u\right \Vert _{C(Q)}=\sup_{(t,x)\in Q}\left \Vert u(t,x)\right \Vert
\text{, \ }\left \Vert u\right \Vert _{C^{\alpha,\beta}(Q)}=\left \Vert
u\right \Vert _{C(Q)}+\sup_{\substack{(t,x),(s,y)\in Q\\t,x)\neq(s,y)}%
}\frac{\left \vert u(t,x)-u(s,y)\right \vert }{\left \vert s-t\right \vert
^{a}+\left \vert x-y\right \vert ^{\beta}},
\]
and
\[%
\begin{array}
[c]{r}%
C_{b}^{1+\alpha,2+\beta}(Q)=\left \{  u:\left \Vert \partial_{t}u\right \Vert
_{C^{\alpha,\beta}(Q)}+\sum \limits_{i=1}^{n}\left \Vert \partial_{x_{i}%
}u\right \Vert _{C^{\alpha,\beta}(Q)}+\sum \limits_{i,j=1}^{n}\left \Vert
\partial_{x_{i}x_{j}}^{2}u\right \Vert _{C^{\alpha,\beta}(Q)}<\infty \right \}  .
\end{array}
\]
If $f,g\in C_{b}^{1+1/2,2+1}([0,T]\times \mathbb{R})$, $\mathbb{\varphi}\in
C_{b}^{2+1}(\mathbb{R})$, and $a(\cdot), b(\cdot)\in C_{b}^{1+1/2}([0,T])$, by
Theorem 6.4.3 in Krylov \cite{Krylov1987} (see also Theorem 4.4 in Appendix C
in Peng \cite{P2010}), then there exists a constant $\alpha \in(0,1)$ such that
for each $\kappa>0$, $u\in C_{b}^{1+\alpha/2,2+\alpha}([0,T-\kappa
]\times \mathbb{R})$.
\end{remark}

\begin{lemma}
\label{lemma3}Assume that $(\phi_{s})_{s\in \lbrack t_{n},t_{n+1}]}$ is a
process in $M_{G}^{2}(t_{n},t_{n+1})$ and $\sup_{s\in \lbrack t_{n,}t_{n+1}%
]}\mathbb{\hat{E}}_{t_{n}}^{x}[\left \vert \phi_{s}\right \vert ^{2}]\leq C$ with 
a constant $C>0$. Then
\begin{equation}
\mathbb{\hat{E}}_{t_{n}}^{x}\left[  \pm \int_{t_{n}}^{t_{n+1}}\int_{t}%
^{t_{n+1}}\phi_{s}dB_{s}d\langle B\rangle_{t}\right]  =0. \label{6.1}%
\end{equation}

\begin{proof}
Let $U_{t}=\int_{t}^{t_{n+1}}\phi_{s}dB_{s}\ $and $U_{t}^{K}=\sum_{i=0}%
^{K-1}\int_{s_{i+1}}^{t_{n+1}}\phi_{s}dB_{s}\mathbf{I}_{[s_{i},s_{i+1})}(t)$,
where $t_{n}=s_{1}<\cdots<s_{K}=t_{n+1}$ is a uniform partition on
$[t_{n,}t_{n+1}]$. From (\ref{isometry}), (\ref{<B><=B}), and Proposition
\ref{PropE}, we can deduce%
\begin{align}
\mathbb{\hat{E}}_{t_{n}}^{x}\left[  \int_{t_{n}}^{t_{n+1}}|U_{t}-U_{t}%
^{K}|^{2}dt\right]   &  \leq \sum \limits_{i=0}^{K-1}\int_{s_{i}}^{s_{i+1}%
}\mathbb{\hat{E}}_{t_{n}}^{x}\left[  \left(  \int_{t}^{s_{i+1}}\phi_{s}%
dB_{s}\right)  ^{2}\right]  dt\nonumber \\
&  \leq \sum \limits_{i=0}^{K-1}\int_{s_{i}}^{s_{i+1}}\mathbb{\hat{E}}_{t_{n}%
}^{x}\left[  \int_{t}^{s_{i+1}}|\phi_{s}|^{2}ds\right]  dt\\
&  \leq C\sum \limits_{i=0}^{K-1}(s_{i+1}-s_{i})^{\frac{3}{2}},\nonumber
\end{align}
which implies that $\mathbb{\hat{E}}_{t_{n}}^{x}[\int_{t_{n}}^{t_{n+1}}%
|U_{t}-U_{t}^{K}|^{2}dt]\rightarrow0$, as $K\rightarrow \infty$. Thus%
\begin{equation}
\int_{t_{n}}^{t_{n+1}}U_{t}d\langle B\rangle_{t}=L_{G}^{2}-\lim_{n\rightarrow
\infty}\int_{t_{n}}^{t_{n+1}}U_{t}^{K}d\langle B\rangle_{t}\text{.}
\label{6.4}%
\end{equation}
Note that $\int_{s_{i+1}}^{t_{n+1}}\phi_{s}dB_{s}$ are independent from
$\Delta \langle B\rangle_{i+1}=\langle B\rangle_{s_{i+1}}-\langle
B\rangle_{s_{i}}$, for $i=0,\ldots,K-1$. This follows that
\begin{equation}
\mathbb{\hat{E}}_{t_{n}}^{x}\left[  \pm \int_{s_{i+1}}^{t_{n+1}}\phi_{s}%
dB_{s}\Delta \langle B\rangle_{i+1}\right]  =\mathbb{\hat{E}}_{t_{n}}%
^{x}\left[  \mathbb{\hat{E}}_{t_{n}}^{x}\left[  \pm y\int_{s_{i+1}}^{t_{n+1}%
}\phi_{s}dB_{s}\right]  _{y=\Delta \langle B\rangle_{i+1}}\right]  =0.
\label{6.2}%
\end{equation}
In view of Proposition \ref{PropE} (ii) and $(\ref{6.2})$, we have
\begin{equation}
\mathbb{\hat{E}}_{t_{n}}^{x}\left[  \pm \int_{t_{n}}^{t_{n+1}}U_{t}^{K}d\langle
B\rangle_{t}\right]  =0\text{.} \label{6.3}%
\end{equation}
From (\ref{6.4}) and (\ref{6.3}), our conclusion follows.
\end{proof}
\end{lemma}

\begin{remark}
\label{remark2}Under the same assumption in Lemma \ref{lemma3}, the conclusion
does no longer hold for the integral $\int_{t_{n}}^{t_{n+1}}\int_{t_{n}}%
^{t}\phi_{s}dB_{s}d\langle B\rangle_{t}$. Here we just give the case $\phi
_{s}=1$ to illustrate the result. Notice that
\[
\mathbb{\hat{E}}_{t_{n}}^{x}\left[  \int_{t_{n}}^{t_{n+1}}\int_{t_{n}}%
^{t}dB_{s}d\langle B\rangle_{t}\right]  =\mathbb{\hat{E}}\left[  \int
_{0}^{\Delta t_{n}}B_{s}d\langle B\rangle_{s}\right]  =\frac{1}{3}%
\mathbb{\hat{E}}[B_{\Delta t_{n}}^{3}]=\frac{1}{3}(\Delta t_{n})^{\frac{3}{2}%
}\mathbb{\hat{E}}[B_{1}^{3}].
\]
Hu \cite{H2012} shows that $\mathbb{\hat{E}}[B_{1}^{3}]>0$, which implies
\[
\mathbb{\hat{E}}_{t_{n}}\left[  \int_{t_{n}}^{t_{n+1}}\int_{t_{n}}^{t}%
dB_{s}d\langle B\rangle_{t}\right]  >0.
\]

\end{remark}

\begin{lemma}
\label{lemma4}Let $R_{y}^{n}$ be the local truncation error derived from
(\ref{RY1})-(\ref{RY2}). Assume $a_{s}$ and $b_{s}$ are bounded processes.
Then for sufficiently small $\Delta t$, it holds that

\begin{description}
\item[(i)] For $\theta_{i}\in \lbrack0,1]$ $\left(  i=1,2\right)  ,$ if
$f,g,\varphi \in C^{\frac{1}{2},1}$, we have%
\begin{equation}
\left \vert R_{y}^{n}\right \vert \leq C(\Delta t)^{\frac{3}{2}}.
\label{truncation RY1}%
\end{equation}

\item[(ii)] In particular, for $\theta_{1}\in \lbrack0,1]$ and $\theta_{2}=0$,
if $f,g,u\in C_{b}^{1,2}$, we have%
\begin{equation}
\left \vert R_{y}^{n}\right \vert \leq C(\Delta t)^{2}. \label{truncation RY2}%
\end{equation}

\end{description}
Here $C>0$ is a constant just depending on $T$, $L$, the upper bounds of $a$
and $b$, the upper bounds of derivatives of $f,$ $g$ and $\varphi$.

\begin{proof}
From the definition of $R_{y_{1}}^{n}$ and $R_{y_{2}}^{n}$, we have
\begin{align}
\left \vert R_{y_{1}}^{n}\right \vert  &  \leq \mathbb{\hat{E}}_{t_{n}}%
^{x}\left[  \int_{t_{n}}^{t_{n+1}}X_{s}^{t_{n}}f_{s}-\theta_{1}X_{t_{n}%
}^{t_{n}}f_{t_{n}}-(1-\theta_{1})X_{t_{n+1}}^{t_{n}}f_{t_{n+1}}ds\right]
\label{7.1}\\
&  \vee \mathbb{\hat{E}}_{t_{n}}^{x}\left[  \int_{t_{n}}^{t_{n+1}}\theta
_{1}X_{t_{n}}^{t_{n}}f_{t_{n}}+(1-\theta_{1})X_{t_{n+1}}^{t_{n}}f_{t_{n+1}%
}-X_{s}^{t_{n}}f_{s}ds\right] \nonumber
\end{align}
and%
\begin{align}
\left \vert R_{y_{2}}^{n}\right \vert  &  \leq \mathbb{\hat{E}}_{t_{n}}%
^{x}\left[  \int_{t_{n}}^{t_{n+1}}X_{s}^{t_{n}}g_{s}-\theta_{2}X_{t_{n}%
}^{t_{n}}g_{t_{n}}-(1-\theta_{2})X_{t_{n+1}}^{t_{n}}g_{t_{n+1}}d\langle
B\rangle_{s}\right] \label{7.3}\\
&  \vee \mathbb{\hat{E}}_{t_{n}}^{x}\left[  \int_{t_{n}}^{t_{n+1}}\theta
_{2}X_{t_{n}}^{t_{n}}g_{t_{n}}+(1-\theta_{2})X_{t_{n+1}}^{t_{n}}g_{t_{n+1}%
}-X_{s}^{t_{n}}g_{s}d\langle B\rangle_{s}\right]  .\nonumber
\end{align}
We just give the estimate of the first part in (\ref{7.1}), and the second
part can be similarly obtained. Note that
\begin{align}
&  \mathbb{\hat{E}}_{t_{n}}^{x}\left[  \int_{t_{n}}^{t_{n+1}}X_{s}^{t_{n}%
}f_{s}-\theta_{1}X_{t_{n}}^{t_{n}}f_{t_{n}}-(1-\theta_{1})X_{t_{n+1}}^{t_{n}%
}f_{t_{n+1}}ds\right] \nonumber \\
&  \leq \theta_{1}\left \{  \int_{t_{n}}^{t_{n+1}}\mathbb{\hat{E}}_{t_{n}}%
^{x}[X_{s}^{t_{n}}(f_{s}-f_{t_{n}})]ds+\int_{t_{n}}^{t_{n+1}}\mathbb{\hat{E}%
}_{t_{n}}^{x}[f_{t_{n}}(X_{s}^{t_{n}}-X_{t_{n}}^{t_{n}})]ds\right \}
\label{7.11}\\
&  \text{ \  \ }+(1-\theta_{1})\left \{  \int_{t_{n}}^{t_{n+1}}\mathbb{\hat{E}%
}_{t_{n}}^{x}[X_{s}^{t_{n}}(f_{s}-f_{t_{n+1}})]ds+\int_{t_{n}}^{t_{n+1}%
}\mathbb{\hat{E}}_{t_{n}}^{x}[f_{t_{n+1}}(X_{s}^{t_{n}}-X_{t_{n+1}}^{t_{n}%
})]ds\right \}  \text{.}\nonumber
\end{align}
(i) If $f,g,\varphi \in C^{\frac{1}{2},1}$, Theorem \ref{Nonlinear Feynman-Kac}
indicates that $Y_{t}=u(t,B_{t})\in C^{\frac{1}{2},1}$. Combining
(\ref{3.2.5}), (\ref{3.1}) and Lemma \ref{lemma1}, it is easy to check that
$\left \vert R_{y_{1}}^{n}\right \vert \leq C(\Delta t)^{\frac{3}{2}}$.
Similarly, we have $\left \vert R_{y_{2}}^{n}\right \vert \leq C(\Delta
t)^{\frac{3}{2}}$, which yields (\ref{truncation RY1}). (ii) Denote
$F(t,B_{t}):=f(t,u(t,B_{t}))$ and $H(t,B_{t}):=g(t,u(t,B_{t}))$. If $f$, $g$,
$u\in C_{b}^{1,2}$, from\ Theorem \ref{Nonlinear Feynman-Kac}, we have $F,H\in
C_{b}^{1,2}$. Applying $G$-It\^{o}'s formula to $X_{t}^{t_{n}}F(t,B_{t})$,
%we have
%\begin{align}
%X_{s}^{t_{n}}F_{s}  &  =F_{t_{n}}+\int_{t_{n}}^{s}(\partial_{t}F_{r}%
%+a_{r}\partial_{x}F_{r})X_{r}^{t_{n}}dr+\int_{t_{n}}^{s}(\frac{1}{2}%
%\partial_{xx}^{2}F_{r}X_{r}^{t_{n}}+b_{r}\partial_{x}F_{r}X_{r}^{t_{n}%
%})d\langle B\rangle_{r}\nonumber \\
%&  \text{\  \  \ }+\int_{t_{n}}^{s}(\partial_{x}F_{r}+b_{r}F_{r})X_{r}^{t_{n}%
%}dB_{r}+\int_{t_{n}}^{s}a_{r}F_{r}X_{r}^{t_{n}}d\tilde{B}_{r}. \label{7.7}%
%\end{align}
%From $\left(  \ref{7.7}\right)  $,
we can deduce
\begin{align}
&  \mathbb{\hat{E}}_{t_{n}}^{x}\left[  \int_{t_{n}}^{t_{n+1}}X_{s}^{t_{n}%
}f_{s}-\theta_{1}X_{t_{n}}^{t_{n}}f_{t_{n}}-(1-\theta_{1})X_{t_{n+1}}^{t_{n}%
}f_{t_{n+1}}ds\right] \nonumber \\
&  \leq \int_{t_{n}}^{t_{n+1}}\mathbb{\hat{E}}_{t_{n}}^{x}\left[  \theta
_{1}\int_{t_{n}}^{s}(\partial_{t}F_{r}+a_{r}\partial_{x}F_{r})X_{r}^{t_{n}%
}dr\right. \nonumber \\
&  \text{ \  \ }-\left.  (1-\theta_{1})\int_{s}^{t_{n+1}}(\partial_{t}%
F_{r}+a_{r}\partial_{x}F_{r})X_{r}^{t_{n}}dr\right]  ds\nonumber \\
&  \text{ \  \ }+\int_{t_{n}}^{t_{n+1}}\mathbb{\hat{E}}_{t_{n}}^{x}\left[
\theta_{1}\int_{t_{n}}^{s}(\frac{1}{2}\partial_{xx}^{2}F_{r}+b_{r}\partial
_{x}F_{r})X_{r}^{t_{n}}d\langle B\rangle_{r}\right. \label{7.2}\\
&  \text{ \  \ }-\left.  (1-\theta_{1})\int_{s}^{t_{n+1}}(\frac{1}{2}%
\partial_{xx}^{2}F_{r}+b_{r}\partial_{x}F_{r})X_{r}^{t_{n}}d\langle
B\rangle_{r}\right]  ds\nonumber \\
&  \text{ \  \ }+\int_{t_{n}}^{t_{n+1}}\mathbb{\hat{E}}_{t_{n}}^{x}\left[
\theta_{1}\int_{t_{n}}^{s}(\partial_{x}F_{r}+b_{r}F_{r})X_{r}^{t_{n}}%
dB_{r}\right. \nonumber \\
&  \text{ \  \ }-\left.  (1-\theta_{1})\int_{s}^{t_{n+1}}(\partial_{x}%
F_{r}+b_{r}F_{r})X_{r}^{t_{n}}dB_{r}\right]  ds\nonumber \\
&  \text{ \  \ }+\int_{t_{n}}^{t_{n+1}}\mathbb{\hat{E}}_{t_{n}}^{x}\left[
\theta_{1}\int_{t_{n}}^{s}a_{r}F_{r}X_{r}^{t_{n}}d\tilde{B}_{r}-(1-\theta
_{1})\int_{s}^{t_{n+1}}a_{r}F_{r}X_{r}^{t_{n}}d\tilde{B}_{r}\right]
ds.\nonumber
\end{align}
Since $F_{r}\in C_{b}^{1,2}$, $X_{r}^{t_{n}}>0$ and $\mathbb{\hat{E}}_{t_{n}%
}^{x}[X_{r}^{t_{n}}]=1$, for any $r\in \lbrack t_{n},t_{n+1}]$, then
\begin{align}
&  \mathbb{\hat{E}}_{t_{n}}^{x}\left[  \int_{t_{n}}^{t_{n+1}}X_{s}^{t_{n}%
}f_{s}-\theta_{1}X_{t_{n}}^{t_{n}}f_{t_{n}}-(1-\theta_{1})X_{t_{n+1}}^{t_{n}%
}f_{t_{n+1}}ds\right] \nonumber \\
&  \leq \theta_{1}\int_{t_{n}}^{t_{n+1}}\mathbb{\hat{E}}_{t_{n}}^{x}\left[
\int_{t_{n}}^{s}\left \vert \partial_{t}F_{r}+a_{r}\partial_{x}F_{r}\right \vert
X_{r}^{t_{n}}dr\right]  ds\nonumber \\
&  \text{ \  \ }+\theta_{1}\int_{t_{n}}^{t_{n+1}}\mathbb{\hat{E}}_{t_{n}}%
^{x}\left[  \int_{t_{n}}^{s}\left \vert \frac{1}{2}\partial_{xx}^{2}F_{r}%
+b_{r}\partial_{x}F_{r}\right \vert X_{r}^{t_{n}}dr\right]  ds\nonumber \\
&  \text{ \  \ }+(1-\theta_{1})\int_{t_{n}}^{t_{n+1}}\mathbb{\hat{E}}_{t_{n}%
}^{x}\left[  \int_{s}^{t_{n+1}}\left \vert \partial_{t}F_{r}+a_{r}\partial
_{x}F_{r}\right \vert X_{r}^{t_{n}}dr\right]  ds\\
&  \text{ \  \ }+(1-\theta_{1})\int_{t_{n}}^{t_{n+1}}\mathbb{\hat{E}}_{t_{n}%
}^{x}\left[  \int_{s}^{t_{n+1}}\left \vert \frac{1}{2}\partial_{xx}^{2}%
F_{r}+b_{r}\partial_{x}F_{r}\right \vert X_{r}^{t_{n}}dr\right]  ds\nonumber \\
&  \leq C(\Delta t_{n})^{2}.\nonumber
\end{align}
This implies that $\left \vert R_{y_{1}}^{n}\right \vert \leq C(\Delta t)^{2}$.
Applying $G$-It\^{o}'s formula to $X_{t}^{t_{n}}H\left(  t,B_{t}\right)  $, if
$\theta_{2}=0$, we can obtain%
\begin{align}
&  \mathbb{\hat{E}}_{t_{n}}^{x}\left[  \int_{t_{n}}^{t_{n+1}}X_{s}^{t_{n}%
}g_{s}-\theta_{2}X_{t_{n}}^{t_{n}}g_{t_{n}}-(1-\theta_{2})X_{t_{n+1}}^{t_{n}%
}g_{t_{n+1}}d\langle B\rangle_{s}\right] \nonumber \\
&  \leq \mathbb{\hat{E}}_{t_{n}}^{x}\left[  -\int_{t_{n}}^{t_{n+1}}\int
_{s}^{t_{n+1}}(\partial_{t}H_{r}+a_{r}\partial_{x}H_{r})X_{r}^{t_{n}%
}drd\langle B\rangle_{s}\right] \nonumber \\
&  \text{ \ }+\mathbb{\hat{E}}_{t_{n}}^{x}\left[  -\int_{t_{n}}^{t_{n+1}}%
\int_{s}^{t_{n+1}}(\frac{1}{2}\partial_{xx}^{2}H_{r}+b_{r}\partial_{x}%
H_{r})X_{r}^{t_{n}}d\langle B\rangle_{r}d\langle B\rangle_{s}\right]
\label{7.4}\\
&  \text{ \ }+\mathbb{\hat{E}}_{t_{n}}^{x}\left[  -\int_{t_{n}}^{t_{n+1}}%
\int_{s}^{t_{n+1}}(\partial_{x}H_{r}+b_{r}H_{r})X_{r}^{t_{n}}dB_{r}d\langle
B\rangle_{s}\right] \nonumber \\
&  \text{ \ }+\mathbb{\hat{E}}_{t_{n}}^{x}\left[  -\int_{t_{n}}^{t_{n+1}}%
\int_{s}^{t_{n+1}}a_{r}F_{r}X_{r}^{t_{n}}d\tilde{B}_{r}d\langle B\rangle
_{s}\right]  .\nonumber
\end{align}
Lemma \ref{lemma3} indicates that
\begin{equation}
\mathbb{\hat{E}}_{t_{n}}^{x}\left[  -\int_{t_{n}}^{t_{n+1}}\int_{s}^{t_{n+1}%
}(\partial_{x}H_{r}+b_{r}H_{r})X_{r}^{t_{n}}dB_{r}d\langle B\rangle
_{s}\right]  =\mathbb{\hat{E}}_{t_{n}}^{x}\left[  -\int_{t_{n}}^{t_{n+1}}%
\int_{s}^{t_{n+1}}a_{r}F_{r}X_{r}^{t_{n}}d\tilde{B}_{r}d\langle B\rangle
_{s}\right]  =0\text{.} \label{7.12}%
\end{equation}
Thus, in the case of $\theta_{2}=0$, together with (\ref{<B><=B}),
(\ref{7.3}), (\ref{7.4}), and (\ref{7.12}), by a simple computation similar to
$\left \vert R_{y_{1}}^{n}\right \vert $, we obtain $\left \vert R_{y_{2}}%
^{n}\right \vert \leq$ $C(\Delta t)^{2}$. But if $\theta_{2}\neq0$, Remark
\ref{remark2} shows that the equation (\ref{7.12}) is not equal to $0$, which
yields $\left \vert R_{y_{2}}^{n}\right \vert \leq$ $C(\Delta t)^{\frac{3}{2}}$
and $\left \vert R_{y}^{n}\right \vert \leq$ $C(\Delta t)^{\frac{3}{2}}$. We
complete\ the proof.
\end{proof}
\end{lemma}

Based on the discussion above, let us show the error estimates of Scheme
\ref{Scheme1}.

\begin{theorem}
\label{theorem2}Let $\left(  Y_{t}\right)  _{0\leq t\leq T}$ be the solution
of the $G$-BSDE (\ref{G-BSDE}) and $Y^{n}(n=N,\ldots,0)$ defined in the Scheme
\ref{Scheme1}. Then for sufficiently small $\Delta t$, it holds that:

\begin{description}
\item[(i)] For $\theta_{i}\in \lbrack0,1]$ $\left(  i=1,2\right)  ,$ if
$f,g,\varphi \in C^{\frac{1}{2},1}$ and $a_{s}=a(s,B_{s})$, $b_{s}=(s,B_{s})$
are bounded stochastic processes in $C^{\frac{1}{2},1}$, we have%
\begin{equation}
\mathbb{\hat{E}[}|Y_{t_{n}}-Y^{n}|^{2}]\leq C\Delta t. \label{theorem4.2.1}%
\end{equation}

\item[(ii)] In particular, for $\theta_{1}\in \lbrack0,1]$ and $\theta_{2}=0,$
if $f,g,u\in C_{b}^{1,2}$ and $a_{s}=a(s)$, $b_{s}=b(s)$ are Lipschitz
continuous functions on $[0,T]$, we have%
\begin{equation}
\mathbb{\hat{E}[}|Y_{t_{n}}-Y^{n}|^{2}]\leq C\left(  \Delta t\right)  ^{2}.
\label{theorem4.2.2}%
\end{equation}

\end{description}

Here $C>0$ is a constant just depending on $T$, $L$, the upper bounds of $a$
and $b$, and the upper bounds of derivatives of $f$, $g$, and $\varphi$.

\begin{proof}
By Lemma \ref{lemma2} and \ref{lemma4}, for sufficiently small $\Delta t$, it
holds that

\begin{description}
\item[(i)] For $\theta_{i}\in \lbrack0,1]$ $\left(  i=1,2\right)  ,$ if
$f,g,\varphi \in C^{\frac{1}{2},1}$, $a_{s}=a(s,B_{s})$ and $b_{s}=b(s,B_{s})$
are bounded stochastic processes in $C^{\frac{1}{2},1}$, we have%
\begin{equation}
\left \vert R_{x}^{n}\right \vert \leq C(\Delta t)^{\frac{3}{2}},\text{
\ }\left \vert R_{y}^{n}\right \vert \leq C(\Delta t)^{\frac{3}{2}}. \label{7.9}%
\end{equation}

\item[(ii)] In particular, for $\theta_{1}\in \lbrack0,1]$ and $\theta_{2}=0,$
if $f,g,u\in C_{b}^{1,2}$, $a_{s}=a(s)$ and $b_{s}=b(s)$ are Lipschitz
continuous functions on $[0,T]$, we have%
\begin{equation}
\left \vert R_{x}^{n}\right \vert \leq C(\Delta t)^{2},\text{ \ }\left \vert
R_{y}^{n}\right \vert \leq C(\Delta t)^{2}. \label{7.10}%
\end{equation}
\end{description}
Using Theorem \ref{theorem1}, the conclusion can be directly obtained.
\end{proof}
\end{theorem}

\begin{remark}
Theorem \ref{theorem2} shows that if $a_{s}=a(s,B_{s})$ and $b_{s}=b(s,B_{s})$
are stochastic processes, Scheme \ref{Scheme1} for solving $G$-BSDE
(\ref{G-BSDE}) admits half-order convergence, for any $\theta_{i}\in
[0,1](i=1,2)$. If $a_{s}$ and $b_{s}$ are deterministic  Lipschitz
continuous functions,
$f,g,\varphi,u$ are smooth enough, Scheme \ref{Scheme1} can reach first-order
convergence rate for $\theta_{1}\in \lbrack0,1]$ and $\theta_{2}=0$.
\textcolor[rgb]{1.00,0.00,0.00}{For classical BSDEs, if the coefficients are smooth enough, $\theta$-schemes
yield a first-order convergence rate in solving $Y$ for any $\theta_{i}\in \lbrack0,1](i=1,2)$, see, for example, \cite[Theorem 3.3]{ZhWP2009} and
\cite[Theorem 1.2]{CR2014}. In comparison, when assuming $\underline{\sigma
}=\bar{\sigma}=\sigma=1$, the $G$-expectation $\mathbb{\hat{E}}$ reduces to
classical expectation $\mathbb{E}$, and\ the $G$-Brownian motion  simplifies
to a standard Brownian motion. This implies that $\langle B,B\rangle_{t}=t$
and
\[
\mathbb{\hat{E}}_{t_{n}}^{x}\left[  \pm \int_{0}^{T}\int_{0}^{t}\phi_{s}dB_{s}d\langle B\rangle_{t}\right]  =\mathbb{E}_{t_{n}}^{x}\left[  \int
_{0}^{T}\int_{0}^{t}\phi_{s}dB_{s}dt\right]  =0\text{,}\] \[
\mathbb{\hat{E}}_{t_{n}}^{x}\left[  \pm \int_{0}^{T}\int_{t}^{T}\phi_{s}dB_{s}d\langle B\rangle_{t}\right]  =\mathbb{E}_{t_{n}}^{x}\left[  \int
_{0}^{T}\int_{t}^{T}\phi_{s}dB_{s}dt\right]  =0\text{,}\]
for any bounded process $\phi_{s}$ in $M_{G}^{2}(0,T)$. In this case, under
the smoothness conditions, Lemma \ref{lemma2}, Lemma \ref{lemma4} and Theorem
\ref{theorem2} show that our scheme solves BSDE with a first-order convergence
rate in $Y$ for any $\theta_{i}\in \lbrack0,1](i=1,2)$, which is consistent
with the classical results.}
\end{remark}

\section{Numerical scheme}

\subsection{\textcolor[rgb]{1.00,0.00,0.00}{Limit theorem approximation method}}%

In this section, we propose an implementable numerical scheme based on Peng's
central limit theorem. For a given $G$-expectation space $(\Omega
,\mathcal{H},\mathbb{\hat{E}})$, let $B_{t}=(B_{t}^{1},\ldots,B_{t}^{d}%
)^{\top}$ be the related $d$-dimensional $G$-Brownian motion and
$\langle \overrightarrow{B}\rangle_{t}:=(b_{11},\ldots,b_{1d},b_{21}%
,\ldots,b_{2d},,\ldots,b_{d1},\ldots,b_{dd})^{\top}$\ be the $d^{2}%
$-dimensional mutual variation vector with $b_{ij}=\langle B^{i},B^{j}%
\rangle_{t}$,\ $i,j=1,2,\ldots,d$.\ In view of \cite[Proposition 3.5.7]%
{P2010}, we know that $(B_{1},\langle \overrightarrow{B}\rangle_{1})$ is a
$\mathbb{R}^{d}\times \mathbb{R}^{d^{2}}$-valued $G$-distributed.

Let $\left \{  \left(  X_{i},Y_{i}\right)  \right \}  _{i=1}^{\infty}$ be an
i.i.d. sequence of $\mathbb{R}^{d}\times \mathbb{R}^{d^{2}}$-valued random
variables on a sublinear expectation space $(\Omega,\mathcal{H},\mathbb{\tilde
{E}})$ satisfying the conditions in Theorem \ref{CLT}. For each $\varphi \in
C(\mathbb{R}^{d}\times \mathbb{R}^{d^{2}})$ satisfying linear growth condition,
Theorem \ref{CLT} yields that
\[
\mathbb{\hat{E}}\left[  \varphi(B_{1},\langle \overrightarrow{B}\rangle
_{1})\right]  =\mathbb{\tilde{E}}\left[  \varphi \left(  \sum_{i=1}^{n}%
\frac{X_{i}}{\sqrt{M}},\sum_{i=1}^{n}\frac{Y_{i}}{M}\right)  \right]
+Err_{\mathbb{\hat{E}},n}%
\]
with the error $Err_{\mathbb{\hat{E}},n}$ converge to 0, as $n\rightarrow
\infty$. For some positive integer $M$, we denote
\begin{equation}
\mathbb{\hat{E}}^{M}\left[  \varphi(B_{1},\langle \overrightarrow{B}\rangle
_{1})\right]  :=\mathbb{\tilde{E}}\left[  \varphi \left(  \sum_{i=1}^{M}%
\frac{X_{i}}{\sqrt{M}},\sum_{i=1}^{M}\frac{Y_{i}}{M}\right)  \right]  .
\label{E}%
\end{equation}
From the independence of $\left \{  \left(  X_{i},Y_{i}\right)  \right \}
_{i=1}^{M}$, we can compute the above expectation in the following way
\[%
\begin{array}
[c]{l}%
\mathbb{\tilde{E}}\left[  \varphi \left(  \sum \limits_{i=1}^{M}\frac{X_{i}%
}{\sqrt{M}},\sum \limits_{i=1}^{M}\frac{Y_{i}}{M}\right)  \right] \\
=\mathbb{\tilde{E}}\left[  \mathbb{\tilde{E}}\left[  \varphi \left(  \frac
{1}{\sqrt{M}}\left(  \sum \limits_{i=1}^{M-1}x_{i}+X_{M}\right)  ,\frac{1}%
{M}\left(  \sum \limits_{i=1}^{M-1}y_{i}+Y_{M}\right)  \right)  \right]
_{(x_{i},y_{i})=(X_{i},Y_{i}),i=1,\cdots,M-1}\right] \\
=\mathbb{\tilde{E}}\left[  \mathbb{\tilde{E}}\left[  \cdots \mathbb{\hat{E}%
}\left[  \varphi \left(  \frac{1}{\sqrt{M}}\left(  \sum \limits_{i=1}^{M-1}%
x_{i}+X_{M}\right)  ,\frac{1}{M}\left(  \sum \limits_{i=1}^{M-1}y_{i}%
+Y_{M}\right)  \right)  \right]  _{\substack{\left(  x_{M-1},y_{M-1}\right)
\\=\left(  X_{M-1},Y_{M-1}\right)  }}\cdots \right]  _{(x_{1},y_{1}%
)=(X_{1},Y_{1})}\right]  .
\end{array}
\label{8.2}%
\]
By choosing appropriate integer $M$ and i.i.d. sequence $\left \{  \left(
X_{i},Y_{i}\right)  \right \}  _{i=1}^{M}$, then we could approximate the
distribution of $(B_{1},\langle \overrightarrow{B}\rangle_{1})$ with the
required accuracy.

\begin{remark}
In the case of $d=1$, for each given $0<\sigma^{2}=-\mathbb{\hat{E}}\left[
-B_{1}^{2}\right]  \leq \mathbb{\hat{E}}\left[  B_{1}^{2}\right]  =1$, we
consider the random variables $(X,Y)$ with following distribution
\[
\left(  X,Y\right)  =\left \{
\begin{array}
[c]{ll}%
(1,\text{ }1) & \text{ }q/2\\
(0,\text{ }0) & 1-q\\
(-1,1) & \text{ }q/2
\end{array}
\right.
\]
where $q\in \{ \sigma^{2},1\}$ is the uncertainty parameter. Define the robust
expectation of $\varphi(X,Y)$
\[
\mathbb{\tilde{E}[}\varphi(X,Y)\mathbb{]}:=\sup_{q\in \{ \sigma^{2},1\}}\left[
\frac{q}{2}\left(  \varphi(1,1)+\varphi(-1,1)\right)  +(1-q)\varphi
(0,0)\right]  ,
\]
for any $\varphi \in C(\mathbb{R}\times \mathbb{R})$. Let $\{(X_{i}%
,Y_{i})\}_{i=1}^{\infty}$ be a sequence of i.i.d. $\mathbb{R}\times \mathbb{R}%
$-valued random variables. Specifically, we have $(X_{1},Y_{1})\overset{d}%
{=}(X,Y)$, $(X_{i+1},Y_{i+1})\overset{d}{=}(X_{i},Y_{i})$, and $(X_{i+1}%
,Y_{i+1})$ is independent from $(X_{1},Y_{1}),\ldots,(X_{i},Y_{i})$ for each
$i\in \mathbb{N}$. It is easy to check that the above sequence satisfies the
conditions of Theorem \ref{CLT} with $\mathbb{\tilde{E}[}\frac{1}{2}aX_{1}%
^{2}+pY_{1}]=\mathbb{\hat{E}[}\frac{1}{2}aB_{1}^{2}+p\langle B\rangle
_{1}]=G\left(  p,a\right)  $, for any $(p,a)\in \mathbb{R}\times \mathbb{R}$.
For the sake of clarity, let us take $M=2$ as an example to give the
implementation of $\mathbb{\hat{E}}^{M}\mathbb{[}\varphi \left(  B_{1},\langle
B\rangle_{1}\right)  ]$
\begin{align*}
&  \mathbb{\hat{E}}^{M}\mathbb{[}\varphi \left(  B_{1},\langle B\rangle
_{1}\right)  ]\\
&  =\mathbb{\tilde{E}}\left[  \mathbb{\tilde{E}[}\tilde{\varphi}(x_{1}%
+X_{2},y_{1}+Y_{2})]_{\left(  x_{1},y_{1}\right)  =\left(  X_{1},Y_{1}\right)
}\right] \\
&  =\sup \limits_{q\in \{ \sigma^{2},1\}}\left \{
\begin{array}
[c]{l}%
\text{ \ }\sup \limits_{q\in \{ \sigma^{2},1\}}\left[  \frac{q}{2}\tilde
{\varphi}(2,2)+(1-q)\tilde{\varphi}(1,1)+\frac{q}{2}\tilde{\varphi
}(0,2)\right]  \frac{q}{2}\\
+\sup \limits_{q\in \{ \sigma^{2},1\}}\left[  \frac{q}{2}\tilde{\varphi
}(1,1)+(1-q)\tilde{\varphi}(0,0)+\frac{q}{2}\tilde{\varphi}(-1,1)\right]
(1-q)\\
+\sup \limits_{q\in \{ \sigma^{2},1\}}\left[  \frac{q}{2}\tilde{\varphi
}(0,2)+(1-q)\tilde{\varphi}(-1,1)+\frac{q}{2}\tilde{\varphi}(-2,2)\right]
\frac{q}{2}%
\end{array}
\right \}  ,
\end{align*}
where $\tilde{\varphi}(x,y)=\varphi(\frac{x}{\sqrt{M}},\frac{y}{M})$.
\end{remark}

\subsection{\textcolor[rgb]{1.00,0.00,0.00}{Fully-discrete scheme}}

Introduce the uniform space partition $\mathcal{D}_{h}=\mathcal{D}_{1,h}%
\times \mathcal{D}_{2,h}\times \cdots \times \mathcal{D}_{d,h}$, where
$\mathcal{D}_{j,h}$ is the partition of the one-dimensional real axis
$\mathbb{R}$
\[
\mathcal{D}_{j,h}=\left \{  \left.  x_{k}^{j}\right \vert x_{k}^{j}=k\Delta
x,\text{ }k=0,\pm1,\pm2,\ldots \right \}  ,
\]
for $j=1,2,\ldots,d$, and $\Delta x$ is a suitable spatial step. Let
$x_{\mathbf{i}}=(x_{i_{1}}^{1},x_{i_{2}}^{2},\ldots,x_{i_{d}}^{d})$\ for
$\mathbf{i}=(i_{1},i_{2},\ldots,i_{d})\in \mathbb{Z}^{d}$ and $\hat
{x}_{\mathbf{i}}=x_{\mathbf{i}}+\Delta B_{n,t_{n+1}}$.
\textcolor[rgb]{1.00,0.00,0.00}{Let $I_{h}u(\hat{x}_{\mathbf{i}})$ be an interpolation approximation of the function $u(x)$ at
the space point $\hat{x}_{\mathbf{i}}$ by using the values of $u(x)$ at a
finite number of the space points
$x_{\mathbf{k}}(\mathbf{k}\in \mathbb{Z}^{d})\in \mathcal{D}_{h}\ $near the space point $\hat{x}_{\mathbf{i}}$. }
We denote $\mathbb{\hat{E}}_{t_{n}}^{x_{\mathbf{i}},M}[\cdot]$ as an
approximation of the conditional mathematical expectation $\mathbb{\hat{E}%
}_{t_{n}}^{x_{\mathbf{i}}}[\cdot]$ using Peng's central limit theorem with
positive integer $M$.

Using the operators $I_{h}$ and $\mathbb{\hat{E}}_{t_{n}}^{x_{\mathbf{i}}%
,M}[\cdot]$, we rewrite the reference equation (\ref{multi-reference}) at the
time-space point $\left(  t_{n},x_{\mathbf{i}}\right)  $ as follows
\begin{align}
Y_{t_{n}}^{x_{\mathbf{i}}}  &  =\mathbb{\hat{E}}_{t_{n}}^{x_{\mathbf{i}}%
,M}\left[  \tilde{X}_{t_{n+1}}^{t_{n}}I_{h}Y_{t_{n+1}}(\hat{x}_{\mathbf{i}%
})+\left(  \theta_{1}f(t_{n},Y_{t_{n}}^{x_{\mathbf{i}}})+(1-\theta_{1}%
)\tilde{X}_{t_{n+1}}^{t_{n}}f\left(  t_{n+1},I_{h}Y_{t_{n+1}}(\hat
{x}_{\mathbf{i}})\right)  \right)  \Delta t_{n}\right. \nonumber \\
&  \text{ \  \ }+\left.  \left(  \theta_{2}g(t_{n},Y_{t_{n}}^{x_{\mathbf{i}}%
})+(1-\theta_{2})\tilde{X}_{t_{n+1}}^{t_{n}}g(t_{n+1},I_{h}Y_{t_{n+1}}(\hat
{x}_{\mathbf{i}})),\Delta \langle B\rangle_{n+1}\right)  \right]
+R^{n}+R_{\mathbb{\hat{E}}}^{n}+R_{I}^{n}, \label{fully_discrete}%
\end{align}
where
\begin{align*}
R_{\mathbb{\hat{E}}}^{n}  &  =\mathbb{\hat{E}}_{t_{n}}^{x_{\mathbf{i}}}\left[
\tilde{X}_{t_{n+1}}^{t_{n}}Y_{t_{n+1}}+\left(  \theta_{1}f(t_{n},Y_{t_{n}%
}^{x_{\mathbf{i}}})+(1-\theta_{1})\tilde{X}_{t_{n+1}}^{t_{n}}f(t_{n+1}%
,Y_{t_{n+1}})\right)  \Delta t_{n}\right. \\
&  \text{ \ }+\left.  \left(  \theta_{2}g(t_{n},Y_{t_{n}}^{x_{\mathbf{i}}%
})+(1-\theta_{2})\tilde{X}_{t_{n+1}}^{t_{n}}g(t_{n+1},Y_{t_{n+1}}%
),\Delta \langle B\rangle_{n+1}\right)  \right] \\
&  -\mathbb{\hat{E}}_{t_{n}}^{x_{\mathbf{i}},M}\left[  \tilde{X}_{t_{n+1}%
}^{t_{n}}Y_{t_{n+1}}+\left(  \theta_{1}f(t_{n},Y_{t_{n}}^{x_{\mathbf{i}}%
})+(1-\theta_{1})\tilde{X}_{t_{n+1}}^{t_{n}}f(t_{n+1},Y_{t_{n+1}})\right)
\Delta t_{n}\right. \\
&  +\left.  \left(  \theta_{2}g(t_{n},Y_{t_{n}}^{x_{\mathbf{i}}}%
)+(1-\theta_{2})\tilde{X}_{t_{n+1}}^{t_{n}}g(t_{n+1},Y_{t_{n+1}}%
),\Delta \langle B\rangle_{n+1}\right)  \right]
\end{align*}
and
\begin{align*}
R_{I}^{n}  &  =\mathbb{\hat{E}}_{t_{n}}^{x_{\mathbf{i}},M}\left[  \tilde
{X}_{t_{n+1}}^{t_{n}}Y_{t_{n+1}}+\left(  \theta_{1}f(t_{n},Y_{t_{n}%
}^{x_{\mathbf{i}}})+(1-\theta_{1})\tilde{X}_{t_{n+1}}^{t_{n}}f(t_{n+1}%
,Y_{t_{n+1}})\right)  \Delta t_{n}\right. \\
&  \text{ \ }+\left.  \left(  \theta_{2}g(t_{n},Y_{t_{n}}^{x_{\mathbf{i}}%
})+(1-\theta_{2})\tilde{X}_{t_{n+1}}^{t_{n}}g(t_{n+1},Y_{t_{n+1}}%
),\Delta \langle B\rangle_{n+1}\right)  \right] \\
&  -\mathbb{\hat{E}}_{t_{n}}^{x_{\mathbf{i}},M}\left[  \tilde{X}_{t_{n+1}%
}^{t_{n}}I_{h}Y_{t_{n+1}}(\hat{x}_{\mathbf{i}})+\left(  \theta_{1}%
f(t_{n},Y_{t_{n}}^{x_{\mathbf{i}}})+(1-\theta_{1})\tilde{X}_{t_{n+1}}^{t_{n}%
}f\left(  t_{n+1},I_{h}Y_{t_{n+1}}(\hat{x}_{\mathbf{i}})\right)  \right)
\Delta t_{n}\right. \\
&  +\left.  \left(  \theta_{2}g(t_{n},Y_{t_{n}}^{x_{\mathbf{i}}}%
)+(1-\theta_{2})\tilde{X}_{t_{n+1}}^{t_{n}}g(t_{n+1},I_{h}Y_{t_{n+1}}(\hat
{x}_{\mathbf{i}})),\Delta \langle B\rangle_{n+1}\right)  \right]  .
\end{align*}

Based on (\ref{fully_discrete}), by omitting numerical errors $R^{n}$,
$R_{\mathbb{\hat{E}}}^{n}$ and $R_{I}^{n}$, we propose the following
fully-discrete $\theta$-scheme for solving $G$-BSDE (\ref{G-BSDE}).

\begin{sch}
\label{fully discrete Scheme multi-dim}
\textcolor[rgb]{1.00,0.00,0.00}{Denote $\hat{Y}_{\mathbf{i}}^{n}$
$(n=N,\ldots,1,0,$ $\mathbf{i}\in \mathbb{Z}^{d}\mathbb{)}$ as a numerical
solution to the analytic solution $Y_{t}$ of the $G$-BSDE (\ref{G-BSDE}) at
the time-space point $(t_{n},x_{\mathbf{i}})$.} Given $\hat{Y}_{\mathbf{i}%
}^{N}$$(\mathbf{i}\in \mathbb{Z}^{d})$, solve $\hat{Y}_{\mathbf{i}}^{n}$,
$\mathbf{i}\in \mathbb{Z}^{d}$, $n=N-1,\ldots,0$, from
\begin{equation}
\left \{
\begin{array}
[c]{l}%
\hat{Y}_{\mathbf{i}}^{n}\text{ }=\mathbb{\hat{E}}_{t_{n}}^{x_{\mathbf{i}}%
,M}\left[  X_{\mathbf{i}}^{n+1}I_{h}\hat{Y}^{n+1}(\hat{x}_{\mathbf{i}%
})+\left(  \theta_{1}f(t_{n},\hat{Y}_{\mathbf{i}}^{n})+(1-\theta
_{1})X_{\mathbf{i}}^{n+1}f(t_{n+1},I_{h}\hat{Y}^{n+1}(\hat{x}_{\mathbf{i}%
}))\right)  \Delta t_{n}\right. \\
\text{ \  \  \  \  \ }+\left.  \left(  \theta_{2}g(t_{n},\hat{Y}_{\mathbf{i}}%
^{n})+(1-\theta_{2})X_{\mathbf{i}}^{n+1}g(t_{n+1},I_{h}\hat{Y}^{n+1}(\hat
{x}_{\mathbf{i}})),\Delta \langle B\rangle_{n+1}\right)  \right]  ,\\
X_{\mathbf{i}}^{n+1}\text{ is the value of }X^{n+1}\text{\ at the time-space
point }(t_{n},x_{\mathbf{i}}),
\end{array}
\right.
\end{equation}
with the deterministic parameters $\theta_{i}\in \left[  0,1\right]  (i=1,2)$.
\end{sch}

\begin{remark}
For the sake of clarity, we take the following special case as an example
\begin{equation}%
\begin{array}
[c]{r}%
\mathbb{\hat{E}}_{t_{n}}^{x_{\mathbf{i}},M}\left[  I_{h}\hat{Y}^{n+1}\left(
B_{t_{n+1}}\right)  \Delta \langle B^{u},B^{v}\rangle_{n+1}\right]
=\mathbb{\tilde{E}}\left[  I_{h}\hat{Y}^{n+1}\left(  x_{\mathbf{i}}%
+\sqrt{\Delta t_{n}}\sum \limits_{m=1}^{M}\frac{X_{m}}{\sqrt{M}}\right)  \Delta
t_{n}\sum \limits_{m=1}^{M}\frac{Y_{m}^{(u-1)d+v}}{M}\right]  ,
\end{array}
\nonumber
\end{equation}
for $u,v=1,\ldots,d$.
\end{remark}

\begin{remark}
The local truncation errors of Scheme \ref{fully discrete Scheme multi-dim}
consist of three terms $R^{n}$, $R_{\mathbb{\hat{E}}}^{n}$, and $R_{I}^{n}$.
$R^{n}$ is the numerical integration approximation error discussed in Lemma
\ref{lemma2} and Lemma \ref{lemma4}. $R_{\mathbb{\hat{E}}}^{n}$ is the
numerical error introduced by approximating conditional expectation.
$R_{I}^{n}$ is the numerical interpolation error. Suppose that $a$, $b$, $f$,
$g$, and $\varphi$ are sufficiently smooth, if the $M$-time summation and the
$r$-degree polynomial interpolation are used, then the following estimates
hold (cf. \cite{AS1972,HL2024,Krylov2018,Song2017})
\[
R_{\mathbb{\hat{E}}}^{n}=\mathcal{O}(1/M)^{1/6}\  \  \  \  \text{and \  \  \ }%
R_{I}^{n}=\mathcal{O}(\Delta x)^{r+1}.
\]
In our numerical experiments, we employ cubic spline interpolation\ with
$r=3$\ and use the sum number $M=20$ in the one-dimensional case, and $M=6$ in
the two-dimensional case.
\end{remark}

\begin{remark}
\textcolor[rgb]{1.00,0.00,0.00}{In practical applications, we are often interested in the values of $\hat
{Y}_{\mathbf{i}}^{0}$ in a bounded domain. In fact,
if $\varphi$ is a bounded function, then
$\mathbb{\hat{E}}_{t_{n}}^{x_{\mathbf{i}},M}[\varphi(B_{t_{n+1}})]=\mathbb{\hat{E}}_{t_{n}
}^{x_{\mathbf{i}},M}\big[  \varphi(x_{\mathbf{i}}+\sqrt{\Delta t_{n}}
B_{1})\mathbf{1}_{\{ \left \vert B_{1}\right \vert \leq D\}}\big]
+Err_{\varphi,D}$, where the truncation error $\left \vert
Err_{\varphi,D}\right \vert \leq C\mathbb{\hat{E}}\big[\mathbf{1}
_{\{ \left \vert B_{1}\right \vert >D\}}\big]  $ is negligible in the case of $D\geq3$
(see Section 4.2 in \cite{PYY2019} for more details). Thus the data of
$\hat{Y}_{\mathbf{i}}^{N}$ within a specified bounded domain is needed under some
required accuracy.
}
\end{remark}

%Notice that the spatial mesh $D_{h}$ is essentially unbounded. However, in
%real computations, we are interested in acquiring certain values of $Y_{t}$ at
%$\left(  t_{n},x_{\mathbf{i}}\right)  $ with $x_{\mathbf{i}}$ belonging to a
%bounded domain, $\mathbf{i}\in \mathbb{Z}^{d}$.

\subsection{Numerical examples}

In this part, some numerical experiments have been carried out to illustrate
the accuracy of our numerical schemes. The first two examples are the case of
one-dimensional $G$-Brownian motion, and the third example is the
two-dimensional case.
%In our examples,
%we let $D=3$ and select a appropriate $M>0$.
%We also apply cubic spline interpolation to compute spatial non-grid points.
Set $x_{0}=0,$ $t_{0}=0$, and $T=1$. In the following tables, $|Y_{0}-Y^{0}|$
denotes the absolute error between the exact and numerical solution for
$Y_{t}$ at $(t_{0},x_{0})$, and CR denotes the convergence rate.

\begin{example}
\label{EX1}Consider the following linear $G$-BSDE
\begin{equation}
\left \{
\begin{array}
[c]{l}%
\displaystyle-dY_{t}=-Y_{t}dt+\frac{1}{2}Y_{t}d\langle B\rangle_{t}%
-Z_{t}dB_{t}-dK_{t},\text{ }0\leq t\leq T,\text{\ }\\
\displaystyle Y_{T}=\exp(T)\sin(B_{T}).
\end{array}
\right.  \label{ex1}%
\end{equation}
The analytic solution of (\ref{ex1}) is $Y_{t}=\exp(t)\sin(B_{t})$,
$Z_{t}=\exp(t)\cos(B_{t})$ and $K_{t}=0$.
%The solution $Y_{t}$ at initial time
%$t_{0}=0$ is $Y_{0}=0$.
Errors and convergence rates for different time partitions and different
$\sigma^{2}$ are listed in Table \ref{Table1}. Figure \ref{Figure1} shows the
convergence rates of the $\theta$-scheme for time step varying from $2^{3}$ to
$2^{7}$ with different parameters $\theta_{i}(i=1,2)$. By contrast, we can
clearly see that the convergence rate of the $\theta$-scheme is highly
correlated with $\theta_{2}$. For $\theta_{1}\in \lbrack0,1]$, the scheme is
half-order convergence with $\theta_{2}\neq0$. When $\theta_{2}=0$, the scheme
can achieve first order, which is consistent with our theoretical analysis.
\begin{table}[ptbh]
\caption{Errors and convergence rates for Example \ref{EX1}.}%
\label{Table1}%
{\footnotesize
\begin{align*}
&
\begin{tabular}
[c]{l|l|c|c|c|c|c|c}\hline
\multicolumn{2}{c|}{$\sigma^{2}=0.25$} & \multicolumn{5}{|c|}{$|Y_{0}-Y^{0}|$}
& CR\\ \hline
\multicolumn{2}{c|}{$N$} & 8 & 16 & 32 & 64 & 128 & \\ \hline
$\theta_{1}=0\ $ & $\theta_{2}=0$ & \multicolumn{1}{|l|}{8.10E-03} &
\multicolumn{1}{|l|}{4.00E-03} & \multicolumn{1}{|l|}{2.00E-03} &
\multicolumn{1}{|l|}{1.00E-03} & \multicolumn{1}{|l|}{4.84E-04} &
\multicolumn{1}{|l}{1.013}\\ \hline
$\theta_{1}=0.5$ & $\theta_{2}=0$ & \multicolumn{1}{|l|}{2.10E-03} &
\multicolumn{1}{|l|}{1.10E-03} & \multicolumn{1}{|l|}{5.17E-04} &
\multicolumn{1}{|l|}{2.90E-04} & \multicolumn{1}{|l|}{1.29E-04} &
\multicolumn{1}{|l}{0.994}\\ \hline
$\theta_{1}=1$ & $\theta_{2}=0$ & \multicolumn{1}{|l|}{3.70E-03} &
\multicolumn{1}{|l|}{1.90E-03} & \multicolumn{1}{|l|}{9.18E-04} &
\multicolumn{1}{|l|}{4.88E-04} & \multicolumn{1}{|l|}{2.28E-04} &
\multicolumn{1}{|l}{1.003}\\ \hline
$\theta_{1}=0\ $ & $\theta_{2}=0.5$ & \multicolumn{1}{|l|}{1.38E-02} &
\multicolumn{1}{|l|}{1.01E-02} & \multicolumn{1}{|l|}{7.20E-03} &
\multicolumn{1}{|l|}{5.20E-03} & \multicolumn{1}{|l|}{3.60E-03} &
\multicolumn{1}{|l}{0.482}\\ \hline
$\theta_{1}=0.5$ & $\theta_{2}=1$ & \multicolumn{1}{|l|}{2.94E-02} &
\multicolumn{1}{|l|}{2.08E-02} & \multicolumn{1}{|l|}{1.47E-02} &
\multicolumn{1}{|l|}{1.04E-02} & \multicolumn{1}{|l|}{7.30E-03} &
\multicolumn{1}{|l}{0.501}\\ \hline
$\theta_{1}=1$ & $\theta_{2}=1$ & \multicolumn{1}{|l|}{3.18E-02} &
\multicolumn{1}{|l|}{2.18E-02} & \multicolumn{1}{|l|}{1.50E-02} &
\multicolumn{1}{|l|}{1.05E-02} & \multicolumn{1}{|l|}{7.40E-03} &
\multicolumn{1}{|l}{0.526}\\ \hline
\end{tabular}
\\
&
\begin{tabular}
[c]{l|l|c|c|c|c|c|c}\hline
\multicolumn{2}{c|}{$\sigma^{2}=0.5$} & \multicolumn{5}{|c|}{$|Y_{0}-Y^{0}|$}
& CR\\ \hline
\multicolumn{2}{c|}{$N$} & 8 & 16 & 32 & 64 & 128 & \\ \hline
$\theta_{1}=0\ $ & $\theta_{2}=0$ & \multicolumn{1}{|l|}{5.30E-03} &
\multicolumn{1}{|l|}{2.60E-03} & \multicolumn{1}{|l|}{1.30E-03} &
\multicolumn{1}{|l|}{6.58E-04} & \multicolumn{1}{|l|}{3.13E-04} &
\multicolumn{1}{|l}{1.012}\\ \hline
$\theta_{1}=0.5$ & $\theta_{2}=0$ & \multicolumn{1}{|l|}{9.98E-04} &
\multicolumn{1}{|l|}{5.43E-04} & \multicolumn{1}{|l|}{2.48E-04} &
\multicolumn{1}{|l|}{1.54E-04} & \multicolumn{1}{|l|}{6.19E-05} &
\multicolumn{1}{|l}{0.984}\\ \hline
$\theta_{1}=1$ & $\theta_{2}=0$ & \multicolumn{1}{|l|}{3.10E-03} &
\multicolumn{1}{|l|}{1.60E-03} & \multicolumn{1}{|l|}{7.61E-04} &
\multicolumn{1}{|l|}{4.11E-04} & \multicolumn{1}{|l|}{1.90E-04} &
\multicolumn{1}{|l}{0.997}\\ \hline
$\theta_{1}=0\ $ & $\theta_{2}=0.5$ & \multicolumn{1}{|l|}{9.90E-03} &
\multicolumn{1}{|l|}{7.20E-03} & \multicolumn{1}{|l|}{5.10E-03} &
\multicolumn{1}{|l|}{3.70E-03} & \multicolumn{1}{|l|}{2.60E-03} &
\multicolumn{1}{|l}{0.486}\\ \hline
$\theta_{1}=0.5$ & $\theta_{2}=1$ & \multicolumn{1}{|l|}{2.07E-02} &
\multicolumn{1}{|l|}{1.47E-02} & \multicolumn{1}{|l|}{1.03E-02} &
\multicolumn{1}{|l|}{7.30E-03} & \multicolumn{1}{|l|}{5.20E-03} &
\multicolumn{1}{|l}{0.501}\\ \hline
$\theta_{1}=1$ & $\theta_{2}=1$ & \multicolumn{1}{|l|}{2.24E-02} &
\multicolumn{1}{|l|}{1.54E-02} & \multicolumn{1}{|l|}{1.06E-02} &
\multicolumn{1}{|l|}{7.40E-03} & \multicolumn{1}{|l|}{5.20E-03} &
\multicolumn{1}{|l}{0.526}\\ \hline
\end{tabular}
\\
&
\begin{tabular}
[c]{l|l|c|c|c|c|c|c}\hline
\multicolumn{2}{c|}{$\sigma^{2}=0.75$} & \multicolumn{5}{|c|}{$|Y_{0}-Y^{0}|$}
& CR\\ \hline
\multicolumn{2}{c|}{$N$} & 8 & 16 & 32 & 64 & 128 & \\ \hline
$\theta_{1}=0\ $ & $\theta_{2}=0$ & \multicolumn{1}{|l|}{2.50E-03} &
\multicolumn{1}{|l|}{1.30E-03} & \multicolumn{1}{|l|}{6.12E-04} &
\multicolumn{1}{|l|}{3.30E-04} & \multicolumn{1}{|l|}{1.52E-04} &
\multicolumn{1}{|l}{1.008}\\ \hline
$\theta_{1}=0.5$ & $\theta_{2}=0$ & \multicolumn{1}{|l|}{2.81E-04} &
\multicolumn{1}{|l|}{1.78E-04} & \multicolumn{1}{|l|}{6.91E-05} &
\multicolumn{1}{|l|}{6.07E-05} & \multicolumn{1}{|l|}{1.73E-05} &
\multicolumn{1}{|l}{0.960}\\ \hline
$\theta_{1}=1$ & $\theta_{2}=0$ & \multicolumn{1}{|l|}{1.90E-03} &
\multicolumn{1}{|l|}{9.77E-04} & \multicolumn{1}{|l|}{4.68E-04} &
\multicolumn{1}{|l|}{2.60E-04} & \multicolumn{1}{|l|}{1.17E-04} &
\multicolumn{1}{|l}{0.989}\\ \hline
$\theta_{1}=0\ $ & $\theta_{2}=0.5$ & \multicolumn{1}{|l|}{5.40E-03} &
\multicolumn{1}{|l|}{3.90E-03} & \multicolumn{1}{|l|}{2.80E-03} &
\multicolumn{1}{|l|}{2.00E-03} & \multicolumn{1}{|l|}{1.40E-03} &
\multicolumn{1}{|l}{0.489}\\ \hline
$\theta_{1}=0.5$ & $\theta_{2}=1$ & \multicolumn{1}{|l|}{1.11E-02} &
\multicolumn{1}{|l|}{7.90E-03} & \multicolumn{1}{|l|}{5.50E-03} &
\multicolumn{1}{|l|}{3.90E-03} & \multicolumn{1}{|l|}{2.80E-03} &
\multicolumn{1}{|l}{0.500}\\ \hline
$\theta_{1}=1$ & $\theta_{2}=1$ & \multicolumn{1}{|l|}{1.19E-02} &
\multicolumn{1}{|l|}{8.20E-03} & \multicolumn{1}{|l|}{5.70E-03} &
\multicolumn{1}{|l|}{4.00E-03} & \multicolumn{1}{|l|}{2.80E-03} &
\multicolumn{1}{|l}{0.524}\\ \hline
\end{tabular}
\end{align*}
}\end{table}

\begin{figure}[ptb]
\centering
%Requires \usepackage{graphicx}
\includegraphics[width=.8\textwidth]{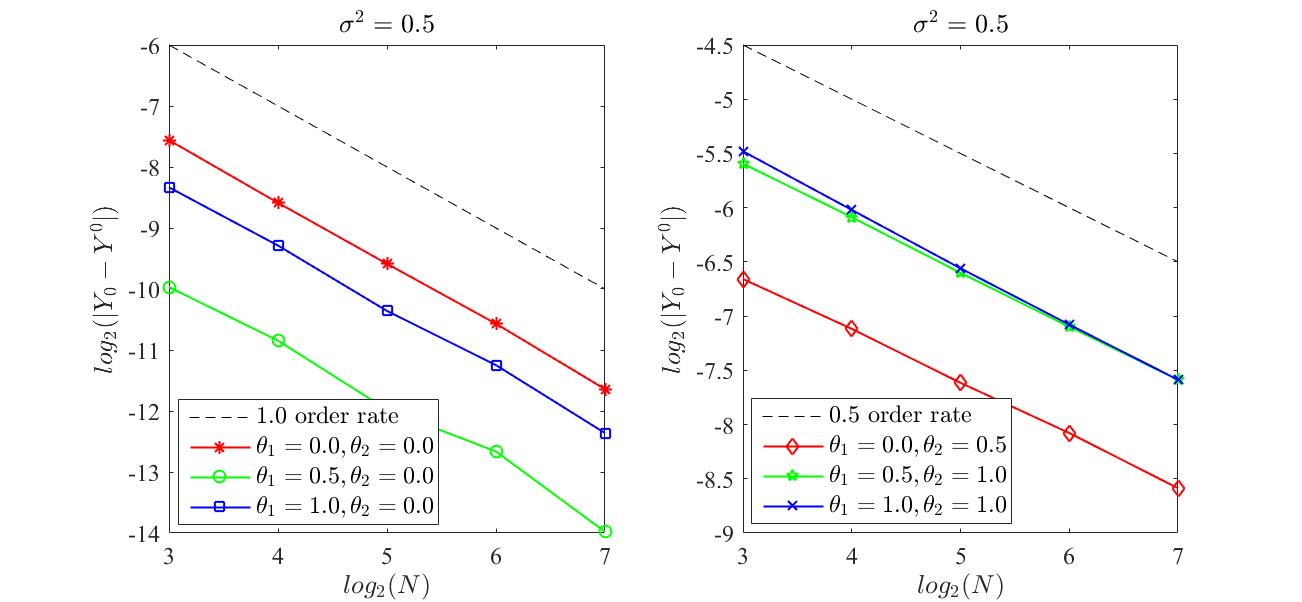}\caption{Convergence rates with
different parameters $\theta_{i}(i=1,2)$. }%
\label{Figure1}%
\end{figure}
\end{example}

\begin{example}
\bigskip \label{EX2}Consider the following nonlinear $G$-BSDE%
\begin{equation}
\left \{
\begin{array}
[c]{l}%
\displaystyle-dY_{t}=\left(  f(t,Y_{t})-Y_{t}\right)  dt+\left(
g(t,Y_{t})+Z_{t}\right)  d\langle B\rangle_{t}-Z_{t}dB_{t}-dK_{t},\text{
}0\leq t\leq T,\\
\displaystyle Y_{T}=\frac{\exp \left(  T+B_{T}\right)  }{1+\exp \left(
T+B_{T}\right)  },
\end{array}
\right.  \label{ex2}%
\end{equation}
where $f(t,y)=J_{\varepsilon}\ast y^{2}$, $g(t,y)=J_{\varepsilon}\ast
(-y^{3}+2.5y^{2}-1.5y)\in C^{\infty}(\mathbb{R})$ with the mollifier$\ $%
\[
J_{\varepsilon}(x)=\frac{k_{0}}{\varepsilon}\exp \left(  -\frac{1}%
{1-|x/\varepsilon|^{2}}\right)  \mathbf{1}_{|x|<\varepsilon}\text{, \ for
given }\varepsilon \in(0,1)\text{ and }k_{0}>0,
\]
and the convolution
\[
J_{\varepsilon}\ast u(x)=\int_{\mathbb{R}}J_{\varepsilon}(x-y)u(y)dy.
\]
The analytic solution of $\left(  \ref{ex2}\right)  $ is $(Y_{t},Z_{t}%
,K_{t})=\left(  \frac{\exp \left(  t+B_{t}\right)  }{1+\exp \left(
t+B_{t}\right)  },\frac{\exp \left(  t+B_{t}\right)  }{\left(  1+\exp \left(
t+B_{t}\right)  \right)  ^{2}},0\right)  $.
%Then $Y_{t}$ at initial time $t_{0}=0$ is $Y_{0}=0.5$.
\begin{table}[ptbh]
\caption{Errors and convergence rates for Example \ref{EX2}.}%
\label{Table2}
{\footnotesize
\begin{align*}
&
\begin{tabular}
[c]{l|l|c|c|c|c|c|c}\hline
\multicolumn{2}{c|}{$\sigma^{2}=0.25$} & \multicolumn{5}{|c|}{$|Y_{0}-Y^{0}|$}
& CR\\ \hline
\multicolumn{2}{c|}{$N$} & 8 & 16 & 32 & 64 & 128 & \\ \hline
$\theta_{1}=0\ $ & $\theta_{2}=0$ & \multicolumn{1}{|l|}{2.52E-02} &
\multicolumn{1}{|l|}{1.29E-02} & \multicolumn{1}{|l|}{6.50E-03} &
\multicolumn{1}{|l|}{3.30E-03} & \multicolumn{1}{|l|}{1.60E-03} &
\multicolumn{1}{|l}{0.988}\\ \hline
$\theta_{1}=0.5$ & $\theta_{2}=0$ & \multicolumn{1}{|l|}{1.30E-02} &
\multicolumn{1}{|l|}{6.60E-03} & \multicolumn{1}{|l|}{3.30E-03} &
\multicolumn{1}{|l|}{1.70E-03} & \multicolumn{1}{|l|}{8.42E-04} &
\multicolumn{1}{|l}{0.988}\\ \hline
$\theta_{1}=0.75$ & $\theta_{2}=0$ & 6.80E-03 & 3.50E-03 & 1.80E-03 &
8.92E-04 & 4.48E-04 & 0.983\\ \hline
$\theta_{1}=0.25$ & $\theta_{2}=0.75$ & \multicolumn{1}{|l|}{5.00E-03} &
\multicolumn{1}{|l|}{3.30E-03} & \multicolumn{1}{|l|}{2.20E-03} &
\multicolumn{1}{|l|}{1.50E-03} & \multicolumn{1}{|l|}{1.10E-03} &
\multicolumn{1}{|l}{0.553}\\ \hline
$\theta_{1}=0.5$ & $\theta_{2}=0.5$ & \multicolumn{1}{|l|}{2.70E-03} &
\multicolumn{1}{|l|}{1.90E-03} & \multicolumn{1}{|l|}{1.30E-03} &
\multicolumn{1}{|l|}{9.58E-04} & \multicolumn{1}{|l|}{6.84E-04} &
\multicolumn{1}{|l}{0.496}\\ \hline
$\theta_{1}=0.15$ & $\theta_{2}=1$ & 5.00E-03 & 3.50E-03 & 2.50E-03 &
1.80E-03 & 1.30E-03 & 0.481\\ \hline
\end{tabular}
\\
&
\begin{tabular}
[c]{l|l|c|c|c|c|c|c}\hline
\multicolumn{2}{c|}{$\sigma^{2}=1$} & \multicolumn{5}{|c|}{$|Y_{0}-Y^{0}|$} &
CR\\ \hline
\multicolumn{2}{c|}{$N$} & 8 & 16 & 32 & 64 & 128 & \\ \hline
$\theta_{1}=0\ $ & $\theta_{2}=0$ & \multicolumn{1}{|l|}{2.52E-02} &
\multicolumn{1}{|l|}{1.28E-02} & \multicolumn{1}{|l|}{6.50E-03} &
\multicolumn{1}{|l|}{3.30E-03} & \multicolumn{1}{|l|}{1.60E-03} &
\multicolumn{1}{|l}{0.988}\\ \hline
$\theta_{1}=0.5$ & $\theta_{2}=0$ & \multicolumn{1}{|l|}{1.30E-02} &
\multicolumn{1}{|l|}{6.60E-03} & \multicolumn{1}{|l|}{3.30E-03} &
\multicolumn{1}{|l|}{1.70E-03} & \multicolumn{1}{|l|}{8.42E-04} &
\multicolumn{1}{|l}{0.988}\\ \hline
$\theta_{1}=0.75$ & $\theta_{2}=0$ & 6.80E-03 & 3.50E-03 & 1.80E-03 &
8.79E-04 & 4.48E-04 & 0.984\\ \hline
$\theta_{1}=0.25$ & $\theta_{2}=0.75$ & \multicolumn{1}{|l|}{1.70E-03} &
\multicolumn{1}{|l|}{8.49E-04} & \multicolumn{1}{|l|}{3.97E-04} &
\multicolumn{1}{|l|}{2.03E-04} & \multicolumn{1}{|l|}{9.74E-05} &
\multicolumn{1}{|l}{1.038}\\ \hline
$\theta_{1}=0.5$ & $\theta_{2}=0.5$ & \multicolumn{1}{|l|}{9.50E-04} &
\multicolumn{1}{|l|}{4.58E-04} & \multicolumn{1}{|l|}{2.03E-04} &
\multicolumn{1}{|l|}{1.06E-04} & \multicolumn{1}{|l|}{4.90E-05} &
\multicolumn{1}{|l}{1.066}\\ \hline
$\theta_{1}=0.15$ & $\theta_{2}=1$ & \multicolumn{1}{|l|}{6.40E-03} &
\multicolumn{1}{|l|}{3.20E-03} & \multicolumn{1}{|l|}{1.50E-03} &
\multicolumn{1}{|l|}{7.77E-04} & \multicolumn{1}{|l|}{3.84E-04} &
1.014\\ \hline
\end{tabular}
\end{align*}
}\end{table}

\begin{figure}[ptbh]
\centering
%Requires \usepackage{graphicx}
\includegraphics[width=.8\textwidth]{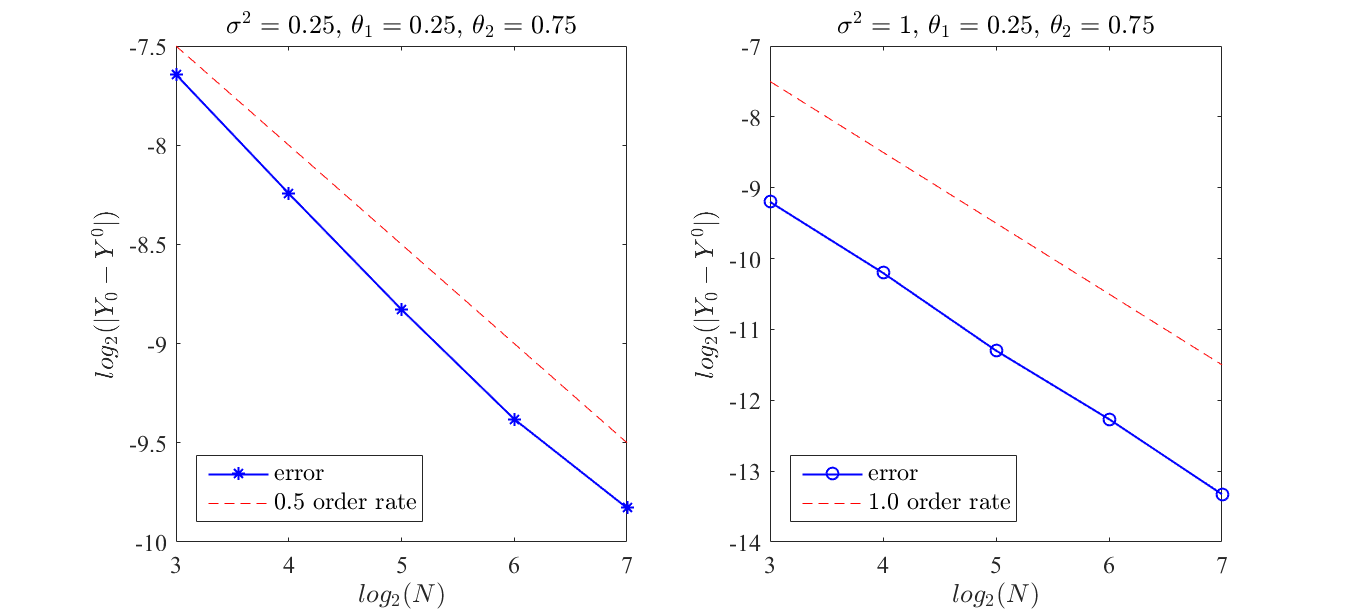}\caption{Convergence rates with
different $\sigma^{2}$. }%
\label{Figure2}%
\end{figure}Errors and convergence rates for different time partitions and
parameters $\theta_{i}(i=1,2)$ are listed in Table \ref{Table2}. Figure
\ref{Table2} presents convergence rates of the $\theta$-scheme under different
uncertainty $\sigma^{2}$. As we can clearly see from Table \ref{Table2},
Scheme \ref{Scheme1} can reach first-order convergence rate when $\theta
_{1}\in \lbrack0,1]$ and $\theta_{2}=0$, and the convergence rate is $1/2$ when
$\theta_{1}\in \lbrack0,1]$ and $\theta_{2}\neq0$. In particular, when
$\underline{\sigma}=\bar{\sigma}=\sigma=1$, for any $\theta_{i}\in \left[
0,1\right]  (i=1,2)$, Scheme \ref{Scheme1} admits a first-order convergence
rate for solving $Y$.
\end{example}

\begin{example}
\label{EX3} let $G(\cdot)$ be a sublinear function defined in \eqref{GA} with
$\Sigma=\{ \lambda Q_{1}+(1-\lambda)Q_{2},\lambda \in \lbrack0,1]\}$, where
\[
Q_{1}=\left(
\begin{array}
[c]{cc}%
1 & 0\\
0 & 1
\end{array}
\right)  \text{, \  \  \ }Q_{2}=\left(
\begin{array}
[c]{cc}%
\frac{1}{2} & 0\\
0 & \frac{1}{4}%
\end{array}
\right)  .
\]
Let $B_{t}=(B_{t}^{1},B_{t}^{2})^{\top}$ be a 2-dimensional $G$-Brownian
motion,
\[
\langle \overrightarrow{B}\rangle_{t}=(\langle B^{1},B^{1}\rangle_{t},\langle
B^{1},B^{2}\rangle_{t},\langle B^{2},B^{1}\rangle_{t},\langle B^{2}%
,B^{2}\rangle_{t})^{\top}%
\]
be the corresponding 4-dimensional mutual variation vector, and $\left(
a_{1},a_{2},a_{3}\right)  :=\left(  1,0,-1\right)  $. We consider
$X=(X^{1},X^{2})^{\top}$ with the distribution
\[
P\{X^{1}=a_{i},X^{2}=a_{j}\}=p_{ij}\text{, for }i,j=1,2,3,
\]
where\ $\tilde{p}:=(p_{ij})_{i,j=1}^{3}\in \{ \tilde{p}_{1},\tilde{p}_{2}\}$
with%
\[
\tilde{p}_{1}=\left(
\begin{tabular}
[c]{ccc}%
$1/4$ & $0$ & $1/4$\\
$0$ & $0$ & $0$\\
$1/4$ & $0$ & $1/4$%
\end{tabular}
\right)  \text{ \ and\  \ }\  \tilde{p}_{2}=\left(
\begin{tabular}
[c]{ccc}%
$1/16$ & $1/8$ & $1/16$\\
$0$ & $1/2$ & $0$\\
$1/16$ & $1/8$ & $1/16$%
\end{tabular}
\right)  .
\]
Let $Y:=(Y^{1},Y^{2},Y^{3},Y^{4})^{\top}$ with $Y^{1}=(X^{1})^{2},$
$Y^{2}=Y^{3}=X^{1}X^{2}$, and $Y^{4}=(X^{2})^{2}$. Define
\[
\mathbb{\tilde{E}[}\varphi(X,Y)\mathbb{]}:=\sup_{\tilde{p}\in \{ \tilde{p}%
_{1},\tilde{p}_{2}\}}\sum_{i,j=1}^{3}\varphi \left(  (a_{i},a_{j})^{\top
},(a_{i}^{2},a_{i}a_{j},a_{i}a_{j},a_{j}^{2})^{\top}\right)  p_{ij}%
,\text{\  \ for }\varphi \in C(\mathbb{R}^{2}\times \mathbb{R}^{4}).
\]
Let $\{(X_{i},Y_{i})\}_{i=1}^{\infty}$ be a sequence of i.i.d. $\mathbb{R}%
^{2}\times \mathbb{R}^{4}$-valued random variables satisfying $(X_{1}%
,Y_{1})\overset{d}{=}(X,Y)$, $(X_{i+1},Y_{i+1})\overset{d}{=}(X_{i},Y_{i})$,
and $(X_{i+1},Y_{i+1})$ is independent from $(X_{1},Y_{1}),\ldots,(X_{i}%
,Y_{i})$ for each $i\in \mathbb{N}$. It can be verified that $\{(X_{i}%
,Y_{i})\}_{i=1}^{\infty}$ satisfies the conditions of Theorem \ref{CLT} with
\[
\mathbb{\tilde{E}}\left[  \frac{1}{2}\langle AX_{1},X_{1}\rangle
+pY_{1}\right]  =\mathbb{\hat{E}}\left[  \frac{1}{2}\langle AB_{1}%
,B_{1}\rangle+p\langle \overrightarrow{B}\rangle_{1}\right]  ,\text{\  \ for
}(p,A)\in \mathbb{R}^{4}\times \mathbb{S}(2).
\]

Now we extend Example \ref{EX1} to the 2-dimensional $G$-Brownian motion case
\begin{equation}
\left \{
\begin{array}
[c]{l}%
\displaystyle-dY_{t}=-Y_{t}dt+\frac{1}{2}Y_{t}d\langle B^{i},B^{j}\rangle
_{t}-Z_{t}dB_{t}-dK_{t},\text{ }0\leq t\leq T,\text{ }i,j=1,2,\text{\ }\\
\displaystyle Y_{T}=\exp(T)\sin(B_{T}^{1}+B_{T}^{2}).
\end{array}
\right.  \label{ex3}%
\end{equation}
The analytic solution of \eqref{ex3} is $Y_{t}=\exp(t)\sin(B_{t}^{1}+B_{t}%
^{2})$, $Z_{t}=\left(  \exp(t)\cos(B_{t}^{1}+B_{t}^{2}),\exp(t)\cos(B_{t}%
^{1}+B_{t}^{2})\right)  $, and $K_{t}=0$. The exact solution is $(Y_{0}%
,Z_{0})=(0,(1,1))$. We list the numerical errors and the convergence rates of
the $\theta$-scheme for time step varying from $2^{3}$ to $2^{7}$ with
different parameters $\theta_{i}(i=1,2)$ in Table \ref{Table5}. It is shown
that our method has high accuracy and aligns well with our theoretical analysis.
\end{example}

\begin{table}[ptbh]
\caption{Errors and convergence rates for Example \ref{EX3}.}%
\label{Table5}%
{\footnotesize
\[%
\begin{tabular}
[c]{l|l|c|c|c|c|c|c}\hline
\multicolumn{2}{c|}{} & \multicolumn{5}{|c|}{$|Y_{0}-Y^{0}|$} & CR\\ \hline
\multicolumn{2}{c|}{$N$} & 8 & 16 & 32 & 64 & 128 & \\ \hline
$\theta_{1}=0\ $ & $\theta_{2}=0$ & \multicolumn{1}{|l|}{7.35E-03} &
\multicolumn{1}{|l|}{3.59E-03} & \multicolumn{1}{|l|}{1.76E-03} &
\multicolumn{1}{|l|}{8.76E-04} & \multicolumn{1}{|l|}{4.36E-04} &
\multicolumn{1}{|l}{1.019}\\ \hline
$\theta_{1}=0.5$ & $\theta_{2}=0$ & \multicolumn{1}{|l|}{1.53E-02} &
\multicolumn{1}{|l|}{7.72E-03} & \multicolumn{1}{|l|}{3.87E-03} &
\multicolumn{1}{|l|}{1.94E-04} & \multicolumn{1}{|l|}{9.68E-04} &
\multicolumn{1}{|l}{0.996}\\ \hline
$\theta_{1}=1$ & $\theta_{2}=0$ & \multicolumn{1}{|l|}{2.84E-02} &
\multicolumn{1}{|l|}{1.47E-02} & \multicolumn{1}{|l|}{7.47E-03} &
\multicolumn{1}{|l|}{3.76E-03} & \multicolumn{1}{|l|}{1.89E-03} &
\multicolumn{1}{|l}{0.979}\\ \hline
$\theta_{1}=0\ $ & $\theta_{2}=0.5$ & \multicolumn{1}{|l|}{3.50E-02} &
\multicolumn{1}{|l|}{2.44E-02} & \multicolumn{1}{|l|}{1.72E-02} &
\multicolumn{1}{|l|}{1.22E-02} & \multicolumn{1}{|l|}{8.60E-03} &
\multicolumn{1}{|l}{0.506}\\ \hline
$\theta_{1}=0.5$ & $\theta_{2}=1$ & \multicolumn{1}{|l|}{6.93E-02} &
\multicolumn{1}{|l|}{4.89E-02} & \multicolumn{1}{|l|}{3.45E-02} &
\multicolumn{1}{|l|}{2.44E-02} & \multicolumn{1}{|l|}{1.72E-02} &
\multicolumn{1}{|l}{0.502}\\ \hline
$\theta_{1}=1$ & $\theta_{2}=1$ & \multicolumn{1}{|l|}{7.21E-02} &
\multicolumn{1}{|l|}{4.99E-02} & \multicolumn{1}{|l|}{3.48E-02} &
\multicolumn{1}{|l|}{2.45E-02} & \multicolumn{1}{|l|}{1.73E-02} &
\multicolumn{1}{|l}{0.515}\\ \hline
\end{tabular}
\]
}\end{table}

\end{document}